\documentclass[a4paper,twoside,12pt]{amsart}

\usepackage[latin1]{inputenc}
\usepackage[english]{babel}

\usepackage{amsxtra}
\usepackage{amssymb}
\usepackage{amsmath,amsthm}
\usepackage[all]{xypic}
\usepackage{epsfig}
\xyoption{dvips}

\theoremstyle{definition}
\newtheorem{ntn}{Notation}[section]
\newtheorem{dfn}[ntn]{Definition}
\newtheorem{rem}[ntn]{Remark}
\newtheorem{exa}[ntn]{Example}
\theoremstyle{plain}
\newtheorem*{que}{Question}
\newtheorem*{thmintr}{Theorem}

\newtheorem{lem}[ntn]{Lemma}

\newtheorem{prp}[ntn]{Proposition}
\newtheorem{thm}[ntn]{Theorem}
\newtheorem{cor}[ntn]{Corollary}
\newtheorem{claim}[ntn]{Claim}

\theoremstyle{remark}

\DeclareMathAlphabet{\mathds}{U}{dsrom}{m}{n}
\DeclareMathAlphabet{\mathsc}{U}{rsfs}{m}{n}

\DeclareMathOperator{\Aut}{Aut}

\DeclareMathOperator{\Gal}{Gal}
\DeclareMathOperator{\Spec}{Spec}

\DeclareMathOperator{\Hg}{Hg}
\DeclareMathOperator{\MT}{MT}

\DeclareMathOperator{\Nm}{Nm}

\DeclareMathOperator{\Pic}{Pic}
\DeclareMathOperator{\End}{End}

\DeclareMathOperator{\CSp}{CSp}

\DeclareMathOperator{\SO}{SO}
\DeclareMathOperator{\GL}{GL}

\newcommand{\Q}{\mathbb{Q}}
\newcommand{\N}{\mathbb{N}}
\newcommand{\Z}{\mathbb{Z}}

\renewcommand{\P}{\mathbb{P}}
\newcommand{\R}{\mathbb{R}}
\newcommand{\A}{\mathbb{A}}
\newcommand{\C}{\mathbb{C}}
\newcommand{\Gm}{\mathbb{G}}
\newcommand{\s}{\mathbb{S}}
\newcommand{\Kg}{\mathbb{K}}

\renewcommand{\O}{\mathcal O}
\renewcommand{\L}{\mathcal L}

\renewcommand{\a}{\alpha}
\renewcommand{\l}{\lambda}

\newcommand{\e}{\epsilon}

\newcommand{\Fk}{\mathcal F}
\newcommand{\Av}{\mathcal A}

\newcommand{\lr}{\rightarrow}

\begin{document}
\title{Complex Multiplication for K3 Surfaces}
\author{Jordan Rizov}
\address{Mathematisch Instituut\\ P.O. Box 80.010\\ 3508 TA Utrecht\\ The Netherlands}
\email{rizov@math.uu.nl}
\begin{abstract}
In this note we prove analogues of the main theorems of complex multiplication for abelian varieties for K3 surfaces. This is done by studying the field of definition of the period morphism for complex K3 surfaces. More precisely we relate the moduli spaces of primitively polarized K3 surfaces with level structures over $\Q$, constructed using algebraic stacks, to the canonical model of the Shimura variety associated to $\SO(2,19)$.
\end{abstract}

\maketitle
\section*{Introduction}

The main theorem for complex multiplication for abelian varieties describes the action of the elements of $\Gal(\Q^{\rm ab}/\Q)$ on the torsion points of an abelian variety $A$ with complex multiplication. In \cite{D-ShV1} Deligne uses this description as a departure point for the
definition of canonical models
of Shimura varieties. In Th\'eor\`eme 4.21 in {\it loc.
cit.} he proves that $\varprojlim \Av_{g,1,n}\otimes \Q$ is the canonical model of $Sh(\CSp_{2g},\mathfrak H^\pm)_\C$, where $\Av_{g,d,n}$ is the moduli stack of $g$-dimensional abelian varieties with a polarization of degree $d^2$ and a Jacobi level $n$-structure. In this note we will prove a similar result for moduli spaces of primitively polarized K3 surfaces.

For a certain class of compact open subgroups $\Kg$ of $\SO(2,19)(\A_f)$, in Section 5 of \cite{Riz-MK3}, we have introduced the notion of a level $\Kg$-structure on a K3 surface using its second \'etale cohomology groups. Further, we have defined moduli spaces $\Fk_{2d,\Kg}$ of
primitively polarized K3 surfaces with a level $\Kg$-structure and we have shown that these are smooth algebraic spaces over $\Spec(\Z[1/N])$ where $N \in \N$
depends on $\Kg$. These moduli spaces are finite unramified covers of the moduli stack $\Fk_{2d}$ of K3 surfaces with a primitive polarization of degree $2d$. 

For each such compact open subgroup $\Kg$, we define below (see Section \ref{ModSpShVAppr}) a period morphism
\begin{displaymath}
j_{d,\Kg, \C} \colon \Fk_{2d,\Kg} \otimes \C \lr Sh_\Kg(\SO(2,19),\Omega^\pm)_\C
\end{displaymath}
where $Sh(\SO(2,19),\Omega^\pm)_\C$ is the Shimura variety associated to $\SO(2,19)$. This period morphism differs slightly from the ones considered for instance Expos\'e XIII in \cite{Ast-K3} and \S 1 in \cite{Fri-Torelli}.
This is due to the fact that we work with moduli spaces of polarized K3 surfaces over $\Q$ and in
general these have more than one connected component. Here we study the field of definition of $j_{2,\Kg,\C}$. Our main result is the following.
\begin{thmintr}
The field of definition of the period morphism $j_{d,\Kg, \C}$ is $\Q$. In other words, $j_{d,\Kg,\C}$ descends to a morphism
\begin{displaymath}
 j_{d,\Kg} \colon \Fk_{2d,\Kg} \otimes \Q \lr Sh_\Kg(\SO(2,19),\Omega^\pm)
\end{displaymath}
where $Sh_\Kg(\SO(2,19),\Omega^\pm)$ is the canonical model of $Sh_\Kg(\SO(2,19),\Omega^\pm)_\C$. 
\end{thmintr}
Just like in the case of abelian varieties it gives a modular interpretation of (an open part of) the canonical model of $Sh(\SO(2,19),\Omega^\pm)_\C$ as a moduli
space.

To prove that the field of definition of $j_{d,\Kg,\C}$ is $\Q$ we need to find `enough' points in
$\Fk_{2d,\Kg} \otimes \C$ and $Sh_\Kg(\SO(2,19),\Omega^\pm)_\C$ for which we can control the action of $\Aut(\C)$. The definition of canonical models of Shimura
varieties suggests a collection of such points, namely the set of special points. Here we will restrict this set a bit by working with points corresponding to
exceptional K3 surfaces. By definition these are K3 surfaces $X$ over $\C$ such that $\text{rk}_\Z\Pic(X)=20$. 

The proof of the preceding theorem splits up in two parts. We first prove a version of the main theorem for complex multiplication (Theorem 11.2 in \cite{Mil-IShV}) for
exceptional K3 surfaces. It gives a relation between the transcendental lattices of an exceptional K3 surface $X$ and of its conjugate $X^\sigma$ by an
automorphism $\sigma \in \Aut(\C)$, fixing the reflex field of $X$. The main tool is a result of Shioda and Inose. Using their construction we reduce the proof of the main theorem of complex multiplication for
exceptional K3 surfaces to a similar statement for abelian surfaces. The latter follows easily from the theorem of Shimura and Taniyama. Next we show that the set of points in $\Fk_{2d,\Kg} \otimes \C$ corresponding to exceptional K3 surfaces with a given reflex field is dense for the Zariski topology. 
The proof of the fact that the field of definition of the period morphism $j_{d,\Kg,\C}$ is $\Q$ is a formal consequence of these two results. 

Using the rationality of the period morphism we show that an analogue of the theorem of Shimura and Taniyama holds for all complex K3 surfaces with complex multiplication. Further, 
we prove that any such K3 surface can be defined over an abelian extension of its Hodge endomorphism field. In this way we complete a theory of complex multiplication for K3 surfaces. 

We should mention that the proofs given here are quite different from the ones given in \cite{Shi-T} in the
case of abelian varieties. Shimura and Taniyama work directly with a given abelian variety $A$ using the geometric interpretation of $H^1_{\rm et}(A,\hat \Z)$.
We obtain our results from the general properties of canonical models of Shimura varieties using the period morphisms. We wonder if one can
give proofs of the statements in this note working directly with an exceptional K3 surface just like in the case of abelian varieties. More
precisely, we wonder if argumentation of the type `a Hodge
cycle is an absolute Hodge cycle' on a K3 surface (see \cite{A-MHC} and \cite{D-HCyc}) could lead to complete proofs of the results in question.
\newline
\newline
{\bf Notations and Conventions}
\newline
\newline
{\bf General.} We write $\hat \Z$ for the profinite completion of $\Z$. We denote by $\A$ the ring of ad\`eles of $\Q$ and by $\A_f = \hat \Z \otimes \Q$ the ring of finite ad\`eles of $\Q$. Similarly, for a number field $E$ we denote by $\A_E$ and $\A_{E,f}$ the ring of ad\`eles and the ring of finite ad\`eles of $E$.

If $A$ is a ring, $A \lr B$ a ring homomorphism then for any $A$-module ($A$-algebra etc.) $V$ we will denote by $V_B$ the $B$-module ($B$-algebra etc.) $V\otimes_A B$.

For a variety $X$ over $\C$ we will denote by $X^{\rm an}$ the associated analytic variety. For an algebraic stack $\mathcal F$ over a scheme $S$ and a morphism of schemes $S' \lr S$ we will denote by $\mathcal F_{S'}$ the product $\mathcal F\times_S S'$ and consider it as an algebraic stack over $S'$.

We will use the notations established in \cite{Riz-MK3}. In particular for a natural number $d$ we write $\Fk_{2d}$ for the Deligne-Mumford stack of K3 spaces with a primitive polarization of degree $2d$. It is a smooth stack over $\Spec(\Z[1/2d])$. See Theorem 4.7 in \cite{Riz-MK3} and \S 1.4.3 in \cite{JR-Thesis}. For $n \in \N$, $n \geq 3$ and a subgroup $\Kg$ of finite index in $\Kg_n$ we denote by $\Fk_{2d,\Kg}$ the smooth algebraic space over $\Spec(\Z[1/N])$ of K3 surfaces with a primitive polarization of degree $2d$ and a level $\Kg$-structure. If $\Kg$ is admissible, then we denote by $\Fk^{\rm full}_{2d,\Kg}$ the moduli space of K3 surfaces with a primitive polarization of degree $2d$ an a full level $\Kg$-structure. For details we refer to Section 6 in \cite{Riz-MK3} and Section 1.5 in \cite{JR-Thesis}.
\newline
\newline 
{\bf Algebraic Groups.} A superscript ${}^0$ usually indicates a connected component for the Zariski topology. For an algebraic group $G$ will denote by $G^0$ the connected component of the identity. We will use the superscript ${}^+$ to denote connected components for other topologies.

For a reductive group $G$ over $\Q$ we denote by $G^{\rm ad}$ the adjoint group of $G$, by $G^{\rm der}$ the derived group of $G$ and by $G^{\rm ab}$ the maximal abelian quotient of $G$. We let $G(\R)_+$ denote the group of elements of $G(\R)$ whose image in $G^{\rm ad}(\R)$ lies in its identity component $G^{\rm ad}(\R)^+$, and we let $G(\Q)_+ = G(\Q) \cap G(\R)_+$. 

Let $V$ be a vector space over $\Q$ and let $G \hookrightarrow \GL(V)$ be an algebraic group over $\Q$. Suppose given a full lattice $L$ in $V$ (i.e., $L\otimes \Q = V$). Then $G(\Z)$ and $G(\hat \Z)$ will denote the abstract groups consisting of the elements in $G(\Q)$ and $G(\A_f)$ preserving the lattices $L$ and $L_{\hat \Z}$ respectively.
\newline
\newline
{\bf Acknowledgments}
\newline
\newline
This note contains the results of Chapter 3 of my Ph.D. thesis \cite{JR-Thesis}. I thank my advisors, Ben Moonen and Frans Oort for their help, their support and for everything I have learned from them. I would like to thank Ben Moonen for pointing out some mistakes in the earlier versions of the text, and Ben Moonen, Gerard van der Geer and Yves Andr\'e for their valuable suggestions. I thank the Dutch Organization for Research N.W.O. for the financial support with which my thesis was done.
\section{Hodge Structures of K3 Surfaces}
\subsection{Hodge Structures}\label{notatch3}
Let $\s = \text{Res}_{\C/\R}(\Gm_{m,\C})$ be the \emph{Deligne torus} over $\R$.
Consider the \emph{weight character} $w \colon \R^\times = \Gm_m(\R) \hookrightarrow \s(\R) = \C^\times$ given by $r \mapsto r^{-1}$. The \emph{norm character} ${\rm Nm} \colon \s \lr \Gm_{m,\R}$ is defined by ${\rm Nm}(z) = z\bar z$. The kernel of ${\rm Nm}$ is the \emph{circle group} ${\mathbb U}_1 = \{z\in\C^* \ |\ |z| = 1\}$.

If a vector space $V$ carries a $\Q$-Hodge structure $h \colon \s \lr \GL(V_\R)$ (shortly $\Q$-HS), then we have a decomposition of $\C$-vector spaces
$$
 V_\C = \bigoplus V^{p,q}
$$
such that $h_\C(z_1,z_2)(v) = z^{-p}_1z_2^{-q}v$ for every $v \in V^{p,q}$. We define $\mu_h$ to be the cocharacter of $\GL(V)_\R$ given by $\mu_h(z) = h_\C(z,1)$.
\begin{exa}\label{muK3}
Let $V$ be a $\Z$-HS of type $\{(1,-1),(0,0),(-1,1)\}$. Then we have a decomposition $V_\C = V^{1,-1} \oplus V^{0,0} \oplus V^{-1,1}$ of $\C$ vector spaces. Then $\mu(z)$ acts on $V^{1,-1}$ as multiplication by $z^{-1}$, on $V^{-1,1}$ as multiplication by $z$ and it acts as the identity on $V^{1,1}$.
\end{exa}
For a $\Q$-HS $V$ we will denote by $\MT(V)$ the \emph{Mumford-Tate group} of $V$. Recall that it is the smallest algebraic subgroup of $\GL(V)$ defined over $\Q$ such that the homomorphism $h$ defining the $\Q$-HS on $V$ factorizes as $h \colon \s \lr \MT(V)_\R \subset \GL(V)_\R$. We will denote by $\Hg(V)$ the \emph{Hodge group} of the $\Q$-HS $V$. By definition it is the smallest algebraic subgroup $\Hg(V) \subset \GL(V)$ such that $h|_{{\mathbb U}_1} \lr \GL(V_\R)$ factorizes through $\Hg(V)_\R$.
\subsection{Quadratic Lattices and Cohomology Groups of K3 Surfaces}\label{QLatK3}
In this section we introduce some notations which will be used in the sequel. Let $U$ be the hyperbolic plane and denote by $E_8$ the \emph{positive} quadratic lattice associated to the Dynkin diagram of type $E_8$ (cf. Chapter V, 1.4 Examples in \cite{S-CA}). 
\begin{ntn}
Denote by $(L_0,\psi)$ the quadratic lattice $U^{\oplus 3} \oplus E_8^{\oplus 2}$. Further, let $(V_0,\psi)$ be the quadratic space $(L_0,\psi) \otimes_\Z \Q$. 
\end{ntn}
\noindent
We have that $L_0$ is a free $\Z$-module of rank 22. The form $\psi_\R$ has signature $(3-,19+)$ on $L_0 \otimes \R$.

Let $\{e_1,f_1\}$ be a basis of the first copy of $U$ in $L_0$ such that
$$
 \psi(e_1,e_1) = \psi(f_1,f_1) = 0\ \ \text{and}\ \ \psi(e_1,f_1) = 1.
$$
For a positive integer $d$ we consider the vector $e_1-df_1$ of $L_0$. It is a primitive vector i.e., the module $L_0/\langle e_1-df_1\rangle$ is free and we have that $\psi(e_1-df_1,e_1-df_1) = -2d$. The orthogonal complement of $e_1-df_1$ in $L_0$ with respect to $\psi$ is $\langle e_1+df_1\rangle \oplus U^{\oplus 2} \oplus E^{\oplus 2}_8$.
\begin{ntn}
Denote the quadratic sublattice $\langle e_1+df_1\rangle \oplus U^{\oplus 2} \oplus E^{\oplus 2}_8$ of $L_0$ by $(L_{2d},\psi_{2d})$. Further, we denote by $(V_{2d},\psi_{2d})$ the quadratic space $(L_{2d},\psi_{2d})\otimes_\Z \Q$. 
\end{ntn}
The signature of the form $\psi_{2d,\R}$ is $(2-,19+)$. We have that $\langle e_1-df_1\rangle \oplus L_{2d}$ is a sublattice of $L_0$ of index $2d$. The inclusion of lattices $i \colon L_{2d} \hookrightarrow L_0$ defines injective homomorphisms of groups
\begin{equation}\label{InjOHom}
i^{\rm ad} \colon \{g \in {\rm O}(V_0)(\Z)\ |\ g(e_1-df_1) = e_1-df_1\}\hookrightarrow {\rm O}(V_{2d})(\Z)
\end{equation}
and
\begin{equation}\label{InjSOHom}
i^{\rm ad} \colon \{g \in \SO(V_0)(\Z)\ |\ g(e_1-df_1) = e_1-df_1\} \hookrightarrow \SO(V_{2d})(\Z).
\end{equation}

Let $X$ be a complex K3 surface. Then one has a non-degenerate bilinear form
$$
 \psi \colon H_B^2(X,\Z)(1) \times H_B^2(X,\Z)(1) \lr \Z
$$
given by
$$
 \psi(x,y) = - {\rm tr} (x \cup y)
$$
where $x \cup y$ is the cup product of $x$ and $y$ and ${\rm tr} \colon H^4_B(X,\Z(2)) \lr \Z$ is the trace map. Let $\L$ be an ample line bundle on $X$ for which $(\L,\L)_X=2d$ and assume that it is primitive. Denote by $P^2_B(X,\Z(1))$ the primitive cohomology group i.e., the orthogonal complement of $c_1(\L)$ in $H^2_B(X,\Z(1))$ with respect to $\psi$. Let $\psi_\L$ be the restriction of $\psi$ to $P^2_B(X,\Z(1))$. Let $\{e_1,f_1\}$ be a basis of the first copy of $U$ in $L_0$. By Proposition 1 in Expos\'e IX in \cite{Ast-K3} one can find an isometry 
\begin{displaymath}
a \colon \bigl(H^2_B(X,\Z(1)),\psi\bigr) \lr L_0
\end{displaymath}
such that $a(c_1(\L)) = e_1 - df_1$. Therefore $a$ induces an isometry 
\begin{displaymath}
a \colon \bigl(P^2_B(X,\Z(1)), \psi_\L\bigr) \lr (L_{2d},\psi_{2d}).
\end{displaymath}
Suppose that $X$ is a K3 surface over a field $k$ of characteristic zero. Then one has a non-degenerate bilinear form
\begin{displaymath}
 \psi_f \colon H_{\rm et}^2(X_{\bar k},\hat \Z)(1) \times H_{\rm et}^2(X_{\bar k},\hat \Z)(1) \lr \hat \Z
\end{displaymath}
given by
$$
 \psi_f (x,y) = -\text{tr}_{\hat \Z}(x\cup y)
$$
where $\text{tr}_{\hat \Z} \colon H^4_{\rm et}(X_{\bar k},\hat \Z)(2) \lr \hat \Z$ is the trace
isomorphism. This is simply Poincar\'e duality for \'etale cohomology (see Corollary 11.2 in \S 11 in \cite{Mil-EC}). One has an isometry
$$
 \bigl(H^2_{\rm et}(X_{\bar k},\hat \Z(1)), \psi_f\bigr) \cong (L_0,\psi)\otimes_\Z \hat \Z.
$$
If $\L$ is a primitive ample line bundle on $X$ of degree $2d$ then one can find a isometry 
\begin{displaymath}
a \colon \bigl(H^2_{\rm et}(X_{\bar k},\hat \Z(1)),\psi_f\bigr) \lr L_0 \otimes_\Z \hat \Z 
\end{displaymath}
such that $a(c_1(\L)) = e_1 - df_1$. Denote again by $P^2_{\rm et}(X_{\bar k},\hat \Z(1))$ the orthogonal complement of $c_1(\L)$ in $H^2_{\rm et}(X_{\bar k},\hat \Z(1))$ with respect to $\psi_f$. Let $\psi_{\L,f}$ be the restriction of $\psi_f$ to $P^2_{\rm et}(X_{\bar k},\hat \Z(1))$. In particular it induces in isometry $\bigl(P^2_{\rm et}(X_{\bar k},\hat \Z(1)), \psi_{\L,f}\bigr) \cong (L_{2d},\psi_{2d})\otimes_\Z \hat \Z$. For a more detailed overview we refer to Section 2 of \cite{Riz-MK3}.
\subsection{Hodge Structures of K3-Type}
In this section we will recall some facts concerning Hodge structures coming form cohomology groups of K3 surfaces. We also set up some notations.

Let $X$ be a K3 surface or an abelian surface over $\C$. Then $H^2_B(X,\Z)$ is a free $\Z$-module of rank either 22 or 6. It carries a $\Z$-HS
$$
 h \colon \s \lr \GL(H^2_B(X,\R))
$$
of weight 2 with Hodge numbers $h^{2,0} = h^{0,2} = 1$ and $h^{1,1} = 20$ or $4$ respectively. Denote by $h_X$ the morphism $h\otimes h_{\Z(1)}$ defining the $\Z$-HS on $H^2_B(X,\Z(1))$.

Assume now that $X$ is a K3 surface and let $\l$ be a primitive quasi-polarization on $X$ (see Section 3.2 in \cite{Riz-MK3} for the definition). Then $P^2_B(X,\Z(1)) = c_1(\l)^\perp \subset H^2_B(X,\Z(1))$ carries a polarized $\Z$-HS. We have that $H^2_B(X,\Q(1)) = c_1(\l) \oplus P^2_B(X,\Q(1))$ as polarized $\Q$-HS. If we consider $\SO\bigl(P^2_B(X,\Q(1))\bigr)$ embedded in to $\SO\bigl(H^2_B(X,\Q(1))\bigr)$ by letting $g \in \SO\bigl(P^2_B(X,\Q(1))\bigr)$ act as the identity on the direct summand $c_1(\l)$ of $H^2_B(X,\Q(1))$, then we have that
$$
 h_X \colon \s \lr \SO\bigl(P^2_B(X,\R(1))\bigr) \hookrightarrow \SO\bigl(H^2_B(X,\R(1))\bigr).
$$
Hence the Mumford-Tate group $\MT\bigl(H^2_B(X,\Q(1))\bigr)$ is contained in $\SO(P^2_B(X,\Q(1)))$. The corresponding cocharacter $\mu_X$, given by $\mu_X(z) = h_{X,\C}(z,1)$, acts as multiplication by $z^{-1}$ on the $\{1,-1\}$ part of $P^2_B(X,\C(1))$, as the identity on the $\{0,0\}$ part and as multiplication by $z$ on the $\{-1,1\}$ part.

Let $X$ be a non-singular projective surface over $\C$ and let $A_X$ be the image 
$$
 A_X := c_1(\Pic(X)) \subset H^2_B(X,\Z(1)).
$$
It is a polarized $\Z$-Hodge substructure of $H^2_B(X,\Z(1))$ of type $\{(0,0)\}$. Denote by $T_X$ its complement 
$$
 T_X := c_1(\Pic(X))^\perp \subset H^2_B(X,\Z(1))
$$
with respect to $\psi_X$ which is the negative of the Poincar\'e duality pairing. The bilinear form $\psi_X$ restricts to a bilinear form on $T_X$. The lattice $T_X$ together with the form $\psi_X|_{T_X}$ is called the \emph{transcendental lattice} of $X$. It carries a polarized $\Z$-HS of type $\{(-1,1),(0,0),(1,-1)\}$ and we have that $H^2_B(X,\Q(1)) = A_{X,\Q} \oplus T_{X,\Q}$ as polarized $\Q$-HS (see \S 1 of \cite{Zar-HGK3} for more details). If we consider $\SO(T_{X,\Q})$ embedded into $\SO\bigl(H^2_B(X,\Q(1))\bigr)$ by acting as the identity on the summand $A_{X,\Q}$ we see that
$$
 h_X \colon \s \lr \SO(T_{X,\R}) \hookrightarrow \SO\bigl(H^2_B(X,\R(1))\bigr).
$$
Hence we have that $\MT\bigl(H^2_B(X,\Q(1))\bigr) \subset \SO(T_{X,\Q})$. The cocharacter $\mu_X$ acts in the way described above.

We fix similar notations for \'etale cohomology groups. Denote by $A_{X,\hat \Z}$ the image $c_1(\Pic(X))_{\hat \Z} \subset H^2_{et}(X,\hat\Z(1))$ and let $A_{X,\A_f}$ be $A_{X,\hat\Z} \otimes \A_f$. We consider the \emph{transcendental lattice} 
$$
 T_{X,\hat \Z} := A_{X,\hat\Z}^\perp \subset H^2_{\rm et}(X,\hat\Z(1))
$$
and let $T_{X,\A_f}$ be $T_{X,\hat\Z}\otimes \A_f$.
Then by the comparison theorem between Betti and \'etale cohomology we have a natural isomorphism $T_{X,\hat\Z} \cong T_X \otimes_\Z \hat\Z$.
\begin{dfn}\label{HSK3TypeDef}
Let $d \in \N$ and consider the vector space $(V_{2d},\psi_{2d})$ defined in Section \ref{QLatK3}. A \emph{polarized $\Q$-HS of K3-type of degree $2d$} is a triple $(V,\psi,h)$ where $h$ is a homomorphism
$$
 h \colon \s \lr \SO(V,\psi)_\R
$$
such that
\begin{enumerate}
\item [(a)] $h$ defines a $\Q$-HS of type $\{(-1,1),(0,0),(1,-1)\}$ on $V$ and $\psi$ is a polarization for $h$,
\item [(b)] The Hodge numbers of $V_\C$ are $h^{-1,1} = h^{1,-1} = 1$ and $h^{0,0} = 19$,
\item [(c)] $(V,\psi)$, as an orthogonal space, is equivalent to $(V_{2d},\psi_{2d})$.
\end{enumerate}
\end{dfn}
Giving such $h$ amounts to giving a $2$-dimensional subspace $V^-_\R$ of $V_\R$ on which $\psi$ is negative definite and an orientation of
$V^-_\R$: For $z = \tau e^{i\theta} \in \s(\R) = \C^\times$, we have that $h(z)$ acts as the identity on the orthogonal complement $V^{-,\perp}_\R$ and as
rotation on angle $2\theta$ on $V^-_\R$ (see Sections 5.4 and 5.5 in \cite{D-K3}).

Note that as above one can see that for a polarized $\Q$-HS of K3-type $h$ the corresponding cocharacter $\mu_h$ acts as in Example \ref{muK3}.
\subsection{Hodge Endomorphism Algebras of K3 Surfaces}
We will give a short exposition of some results of Zarhin (\cite{Zar-HGK3}) on Hodge groups of K3 surfaces. Let $X$ be a complex K3 surface and denote by $\Hg(X)$ and $\MT(X)$ the Hodge and the Mumford-Tate
groups (which in our case are the same) of $X$ associated to the $\Z$-HS on $H^2_B(X,\Z(1))$. We know that the homomorphism
$$
 h_X \colon \s \lr \SO\bigl(H^2_B(X,\R(1))\bigr)
$$
defining the $\Z$-HS on $H^2_B(X,\Z(1))$ factorizes thorough $\SO(T_{X,\R})$ and hence we have that $\Hg(X) \subset \SO(T_{X,\Q})$. We see from Theorem 1.4.1 in \cite{Zar-HGK3} that the vector space $T_{X,\Q}$ is a simple $\Hg(X)$-module.
\begin{dfn}
For a K3 surface $X$ we define the \emph{Hodge endomorphism algebra} of $X$ to be
\begin{displaymath}
 E_X := \End_{\Hg(X)}(T_{X,\Q}) = \End_{\rm HS}(T_{X, \Q}).
\end{displaymath}
\end{dfn}
The Hodge endomorphism algebra of $X$ is an analogue of the endomorphism algebra of an abelian surface.
\begin{thm}
The algebra $E_X$ is a field which is either a totally real field or an imaginary quadratic extension of a totally real field i.e., a CM-field.
\end{thm}
\begin{proof}
For the first part of the theorem we refer to Theorem 1.6 in \cite{Zar-HGK3}. It actually follows from the fact that one has a natural embedding
\begin{displaymath}
\epsilon_X \colon E_X \hookrightarrow \End_\C H^{2,0}(X) \cong \C.
\end{displaymath}
For the second part see Theorem 1.5.1 in \cite{Zar-HGK3}.
\end{proof}

If $X$ is a K3 surface with a Hodge endomorphism algebra $E_X$ which is a CM-field then $E_X$ embedded via $\epsilon_X$ in $\C$ as above is the \emph{reflex field} of the Mumford-Tate torus $\MT(X)$. We will call $\e_X(E_X) \subset \C$ the \emph{reflex field} of $X$.
\begin{dfn}
We say that a K3 surface $X$ has \emph{complex multiplication}, shortly \emph{CM}, by $E_X$ if $X$ has a Hodge endomorphism algebra $E_X$ which is a CM-field.
\end{dfn}
Recall that a complex K3 surface $X$ is called \emph{exceptional} if $\text{rk}_\Z\Pic(X) = 20$. If $X$ is an exceptional K3 surface, then $T_{X, \Q}$ is a 2-dimensional $\Q$-vector space. From \S 2, 2.1 - 2.3 in \cite{Zar-HGK3} we see that the Hodge endomorphism field $E_X$ is a quadratic imaginary field. Further, we have that $\Hg(X)$ (respectively $\MT(X)$) is a torus and therefore an exceptional K3 surface has CM by a quadratic imaginary field.
\begin{dfn}
If $X$ is an exceptional K3 surface over $\C$ we will say that it is of \emph{CM-type} $(E_X,\e_X)$ if $E_X$ is its Hodge endomorphism field and $\e_X \colon E_X \lr \C$ is the natural embedding given by the action of $E_X$ on the space of holomorphic forms on $X$.
\end{dfn}
\section{The Shimura Variety $Sh(\SO(2,19),\Omega^\pm)$}
Over $\C$ the geometry of moduli spaces of primitively polarized K3 surfaces is connected to the geometry of the Shimura variety associated to $\SO(2,19)$. This is achieved via periods. Here we will use this relation to prove the main theorems of complex multiplication for K3 surfaces. In Sections \ref{SOGroups} - \ref{ModSpShVAppr} we will give the needed preliminaries. More precisely, in Section \ref{MIntP1Sect} we describe a modular interpretation of the points in $Sh(\SO(2,19),\Omega^\pm)_\C(\C)$ in terms of periods of K3 surfaces. In the following section we define a period morphism $j_{d,\Kg,\C}$ which is a slight modification of the period morphisms used in \cite{PSh-Sh}, \cite{Fri-Torelli}, \cite{Ast-K3} and others.
\subsection{Special Orthogonal Groups}\label{SOGroups}
Let $V$ be a 21 dimensional $\Q$ vector space and let $\psi$ be a non-degenerate form on it of signature $(2-,19+)$. Then $\psi$ is equivalent to the form 
$$
Q_d \colon -x_1^2 - x_2^2 + x_3^2 + \dots + x^2_{20} + dx_{21}^2
$$
for some square free integer $d$. In general, if $d > 1$ the two forms $Q_d$ and $Q_1$ need not be equivalent over $\Q$. But the forms $dQ_1$ and $Q_d$ are equivalent over $\Q$. One sees this using Theorem and the Corollary in Chapter IV, Section 3.3 in \cite{S-CA}.
Therefore we have the following result.
\begin{lem}
The groups $\SO(V,Q_1)$ and $\SO(V,Q_d)$ are isomorphic over $\Q$.
\end{lem} 
\begin{ntn}\label{GrNotat}
From now on we will denote by $G$ or by $\SO(2,19)$ the algebraic group $\SO(V,\psi)$ over $\Q$.
\end{ntn}
\subsection{The Shimura Datum $(G,\Omega^\pm)$}\label{OrthShimDatum}
Let $V$ be a 21 dimensional vector space over $\Q$ and let $\psi$ be a non-degenerate bilinear form on $V$ of signature $(2-,19+)$. Define $\Omega^\pm$ to be 
the collection of all $\Q$-HS $h\colon \s \lr G_\R = \SO(V_\R)$ for which $\pm\psi$ is a polarization and the Hodge numbers are $h^{-1,1} = h^{1,-1} = 1$ and 
$h^{0,0} = 19$. Let $h_0 \colon \s \lr G_\R$ be an element of $\Omega^\pm$. Then $\Omega^\pm$ is equal to the $G(\R)$-conjugacy class of the homomorphism $h_0$. 

For an element $h\in\Omega^\pm$ let $F_h$ be the associated Hodge
filtration. Then the map $h \mapsto F^1_h$ identifies $\Omega^\pm$ with the space
$$
\{ \omega \in \P(V\otimes \C) |\ \psi(\omega,\omega) = 0\ \text{and}\ \psi(\omega,\bar \omega) > 0 \}.
$$
This gives $\Omega^\pm$ a complex structure for which the Hodge filtration $F_h$ varies
holomorphically with $h$. We also see that $\Omega^\pm$ consists of two connected components $\Omega^+$ and $\Omega^-$ corresponding to the two possible 
orientations one can give to the space $V^-_\R$ corresponding to a morphism $h$. Further $\Omega^\pm$, being the $G(\R)$-conjugacy class of a homomorphism $h_0 \colon \s \lr G_\R$ as above, can be identified with the space 
$$
\SO(2,19)(\R)/\bigl(\SO(2)(\R) \times \SO(19)(\R)\bigr).
$$
We choose $\Omega^+$ to be connected component corresponding to 
$$
 \SO(2,19)(\R)^+/\bigr(\SO(2)(\R) \times \SO(19)(\R)\bigr),
$$
where $\SO(2,19)(\R)^+$ is the connected component of $\SO(2,19)(\R)$ containing the identity. This choice is non-canonical as it depends on the choice of $h_0$.

The pair $(G,\Omega^\pm)$ is a Shimura datum with reflex field $\Q$. The last claim follows from Proposition 3.8 of \cite{D-ShV1} and Appendix 1, Lemma in \cite{A-HV}.
\subsection{Modular Interpretation: Part I}\label{MIntP1Sect}
Let $d \in \N$ and consider the quadratic space $(V_{2d},\psi_{2d})$ (see Section \ref{QLatK3}). Following Notation \ref{GrNotat} we denote  by $G$ the group $\SO(V_{2d},\psi_{2d})$. Let $\Kg \subset G(\A_f)$ be an compact open subgroup and consider the variety
$$
 Sh_\Kg(G,\Omega^\pm)_\C = G(\Q)\backslash \Omega^\pm \times G(\A_f)/\Kg
$$
where $q(x,a)k = (qx,qak)$ for $q \in G(\Q), x \in \Omega^\pm, a \in G(\A_f)$ and $k \in \Kg$.

Define $\mathcal H_\Kg$ to be the set of 4-tuples
$$ 
 \bigl((W,h),s,\a\Kg\bigr)
$$
where:
\begin{enumerate}
\item [(i)]   $\bigl((W,h),s\bigr)$ is a polarized $\Q$-HS of K3 type (see Definition \ref{HSK3TypeDef}),
\item [(ii)] $\a\Kg$ is the $\Kg$-orbit of an $\A_f$-linear isomorphism
$$
 \a \colon V_{2d}\otimes \A_f \lr W \otimes \A_f
$$
such that $\psi_{2d}(v_1,v_2) = m_s \cdot s(\a(v_1),\a(v_2))$ for all $v_1,v_2 \in V_{2d}\otimes \A_f$, where $m_s \in \Q_{+}$.
\end{enumerate}

An isomorphism between $\bigl((W,h),s,\a\Kg\bigr)$ and $\bigl((W',h'),s',\a'\Kg\bigr)$ in $\mathcal H_\Kg$ is an isomorphism $b \colon (W,h) \lr (W',h')$ of $\Q$-HS such that there exists $c \in \Q^\times$ for which $s = c\cdot s'\circ (b \times b)$ and $b\circ \a \equiv \a' \pmod \Kg$. From this we see that $c = m_s/m_{s'} \in \Q^\times/\Q^{\times 2}$.

Let $\bigl((W,h),s,\a\Kg\bigr)$ be an element of $\mathcal H_\Kg$. From (ii) and the fact that the signature of $s$ is $(2-,19+)$ we conclude by the Hasse principle that there is an isomorphism $a \colon W \lr V_{2d}$ with $s = m^{-1}_s\cdot \psi_{2d}\circ (a \times a)$. Further, the homomorphism $a\cdot h \colon \s \lr G_\R$ defined by $z \mapsto a\circ h(z) \circ a^{-1}$ belongs to $\Omega^\pm$ and the
composition
$$
\xymatrix{
 V_{2d, \A_f} \ar[r]^\a & W_{\A_f} \ar[r]^a  & V_{2d, \A_f}
}
$$
is an element of $G(\A_f)$. Indeed, we have that 
\begin{displaymath} 
\psi_{2d}\circ(a\times a)\circ(\a \times \a) = m_s\cdot s \circ (\a \times \a) = m_s m_s^{-1}\cdot \psi_{2d} = \psi_{2d}.
\end{displaymath}
Another isomorphism $a' \colon W \lr V_{2d}$, for which $s = m_s^{-1} \cdot \psi_{2d}\circ (a'\times a')$, differs from $a$ by an element $q\in G(\Q)$, say $a' = q\circ a$. Hence,
replacing $a$ by $a'$ will change $[a\cdot h,a\circ \a]_\Kg$ with $[qa\cdot h,qa\circ \a]_\Kg$. Similarly, replacing $\a$ by $\a k$ for some $k \in \Kg$ will replace
$[a\cdot h,a\circ \a]_\Kg$ with $[a\cdot h,a\circ \a k]_\Kg$. Therefore one has a well defined map
$$
 \mathcal H_\Kg \lr G(\Q)\backslash \Omega^\pm \times G(\A_f)/\Kg
$$
given by
\begin{equation}\label{ratpmap}
 \bigl((W,h),s,\a\Kg\bigr) \mapsto [a\cdot h,a\circ \a]_\Kg
\end{equation}
where $[a\cdot h,a\circ \a]_\Kg$ denotes the class of $(a\cdot h,a\circ \a)$.
\begin{prp}\label{Modinterpr1}
The map defined by \eqref{ratpmap} gives a bijection between $\mathcal H_\Kg/\{{\rm isom.}\}$ and $Sh_\Kg(G,\Omega^\pm)_\C(\C)$.
\end{prp}  
\begin{proof}
Suppose that $b \colon (W,h) \lr (W',h')$ is an isomorphism of $\Q$-HS giving an isomorphism of the triples $\bigl((W,h),s,\a\Kg\bigr)$ and $\bigl((W',h'),s',\a'\Kg\bigr)$ in
$\mathcal H_\Kg$. Then we have that $s = m_s/m_{s'} \cdot s'\circ (b\times b)$. Choose an isomorphism $a' \colon W' \lr V_{2d}$ such that $\psi_{2d} = m_{s'}\cdot s'\circ (a'\times a')$. Then for the isomorphism $a \colon W \lr V_{2d}$ defined by $a = a'\circ b$ we see that $\psi_{2d} = m_s \cdot s \circ (a\times a)$. Hence we have that $(a\cdot h,a\circ \a) = (a'\cdot h',a'\circ \a'k)$ where $b\circ \a = \a'k$.

Assume that $\bigl((W,h),s,\a\Kg\bigr)$ and $\bigl((W',h'),s',\a'\Kg\bigr)$ are mapped to the same point in $Sh_\Kg(G,\Omega^\pm)(\C)$. Choose two isomorphisms $a \colon W \lr V_{2d}$
such that $\psi_{2d} = m_s \cdot s\circ(a\times a)$ and $\a' \colon W' \lr V_{2d}$ for which $\psi_{2d} = m_{s'}\cdot s'\circ(a'\times a')$. We know that 
$$
(a\cdot h,a\circ \a) = (qa'\cdot h',qa'\circ \a'k)
$$ 
for some $q\in G(\Q)$ and $k \in \Kg$. After replacing $a'$ by $qa'$ we may suppose that $(a\cdot h,a\circ \a) = (a'\cdot h',a'\circ \a k)$. Then $b = a'\circ a^{-1}$ is an isomorphism of the triples $\bigl((W,h),s,\a\Kg\bigr)$
and $\bigl((W',h'),s',\a'\Kg\bigr)$. This shows that the map is injective. The surjectivity follows easily as any element $[h,g]_\Kg$ is the image of $\bigl((V_{2d},h),\psi_{2d},g\Kg\bigr)$.
\end{proof}
\begin{rem}
With this modular interpretation of $Sh_\Kg(G,\Omega^\pm)_\C(\C)$ we may think of it as the set parameterizing `isogeny' classes of polarized K3 surfaces with certain
level structure up to an isomorphism. This is similar to the case of abelian varieties. See Section 4.11 in \cite{D-ShV1}.
\end{rem}
\subsection{Modular Interpretation: Part II}\label{ModSpShVAppr}
Let $d$ be a natural number and let $\Kg \subset \SO(V_{2d})(\hat \Z)$ be a subgroup of finite index which is contained in some $\Kg_n$ for some $n \geq 3$. Recall that in \S 6 in \cite{Riz-MK3} (or \S 1.5 in \cite{JR-Thesis}) for every such subgroups $\Kg$ we constructed a moduli space $\Fk_{2d,\Kg}$ of K3 surfaces with a primitive polarization of degree $2d$ and a level $\Kg$-structure. Our goal in this section is to define a morphism
\begin{displaymath}
 j_{d,\Kg,\C} \colon \Fk_{2d,\Kg,\C} \lr  Sh_\Kg(G,\Omega^\pm)_\C
\end{displaymath}
mapping every primitively polarized complex K3 surface its periods (cf. Step 1 of the proof of Proposition \ref{PerMorphism} and Definition \ref{PeriodsDef} below). We will use the notations established in Section \ref{QLatK3}.
\begin{prp}\label{PerMorphism}
For a natural number $d$ and a group $\Kg$ as above one has an \'etale morphism of algebraic spaces
\begin{displaymath}
 j_{d,\Kg,\C} \colon \Fk_{2d,\Kg, \C} \lr Sh_\Kg(G,\Omega^\pm)_\C.
\end{displaymath}
\end{prp} 
\begin{proof} We will divide the proof into several steps.
\newline
{\bf Step 1:} We begin with a naive pointwise definition. Let $(X,\l,\a)$ be a complex K3 surface with a primitive polarization of degree $2d$ and a level $\Kg$-structure $\a$ (see Section 5 of \cite{Riz-MK3} for the definition). Let $\tilde \a \colon L_{2d,\hat \Z} \lr P^2_{\rm et}(X,\hat \Z(1))$ be a representative of the class $\a$. Choose an isometry $a \colon H^2_B(X,\Z(1)) \lr L_0$ such that $a(c_1(\l)) = e_1 - df_1$ and $a\circ \tilde \a \colon L_{2d,\hat \Z} \lr L_{2d, \hat \Z}$ is an element in $\SO(V_{2d})(\hat \Z)$. Let $h_X \colon \s \lr \SO(P^2_B(X,\R(1)))$ be the morphism defining the polarized $\Z$-HS on $P^2_B(X,\Z(1))$.
\begin{claim}\label{CorrectMap}
The class $[a\circ h_X \circ a^{-1}, a\circ \tilde \a]_\Kg$ of the pair $(a\circ h_X \circ a^{-1}, a\circ \tilde \a)$ in $Sh_\Kg(G,\Omega^\pm)_\C(\C)$ is independent of a choice of the marking $a$ and the lifting $\tilde \a$ of $\a$.
\end{claim}
\begin{proof}
Indeed, any representative of the class of $\a$ is of the form $\tilde \a \circ \kappa$ for some $\kappa \in \Kg$ and any isometry $a' \colon H^2_B(X,\Z(1)) \lr L_0$ such that $a(c_1(\l)) = e_1 - df_1$ and $a'\circ \tilde \a \circ \kappa \colon L_{2d,\hat \Z} \lr L_{2d, \hat \Z}$ is equal to $g \circ a$ for some $g \in {\rm O}(V_0)(\Z)$ with $g(e_1 - df_1) = e_1 - df_1$ and such that $g \in \SO(V_{2d})(\Z)$. Hence we have that the new data produce a pair 
\begin{displaymath}
 (g\circ a \circ h_X \circ a^{-1} \circ g^{-1}, g\circ a \circ \tilde \a \circ \kappa) = 
 \bigl(g \cdot (a \circ h_X \circ a^{-1}), g \circ a\circ \tilde \a \circ \kappa\bigr)
\end{displaymath}
whose class in $Sh_\Kg(G,\Omega^\pm)_\C(\C)$ is exactly $[a\circ h_X \circ a^{-1}, a\circ \tilde \a]_\Kg$.
\end{proof}
We will use this pointwise construction to define an algebraic morphism as claimed in the proposition.
\newline
\newline
{\bf Step 2:}  Let $U \lr \Fk_{2d,\Kg, \C}$ be a (smooth) atlas of $\Fk_{2d,\Kg, \C}$ such that the pull-back of the universal family over $\Fk_{2d,\Kg, \C}$ to $U$ is a K3 scheme. Let $V$ be a connected component of $U$ and let $(\pi \colon X \lr V, \l,\a)$ be the pull-back of the universal family to $V$. Define a map
$$
 j_{d,\Kg,V} \colon V^{\rm an} \lr Sh_\Kg(G,\Omega^\pm)_\C
$$
by sending a point $s \in V^{\rm an}$ to the point associated to $(X_s,\l_s,\a_s)$ in Step 1. We will show that it is an algebraic morphism.
\newline
\newline
{\bf Step 3:} We will show that $j_{d,\Kg,V}$ is holomorphic and a local isomorphism. According to Lemma 5.13 of \cite{Mil-IShV} the decomposition of $Sh_\Kg(G,\Omega^\pm)_\C$ into connected components is given in the following way:
\begin{displaymath}
Sh_{\Kg}(G,\Omega^\pm)_\C = \coprod_{[g] \in \mathcal C} \Gamma_{[g]} \backslash \Omega^+,
\end{displaymath}
where $\mathcal C := G(\Q)_+\backslash G(\A_f)/\Kg$ and $\Gamma_{[g]}= G(\Q)_+ \cap g \Kg g^{-1}$ for some representative $g$ of $[g] \in \mathcal C$. We will first show that $j_{d,\Kg,V}$ maps $V^{\rm an}$ into one connected component. 

Suppose given two points $s_1$ and $s_2$ of $V^{\rm an}$. One can find an isomorphism 
$$
\delta_\pi \colon \pi_1(V^{\rm an},s_1) \cong \pi_1(V^{\rm an},s_2)
$$ 
and an isometry 
$$
\delta_B \colon H^2_B(X_{s_1}, \Z(1)) \lr H^2_B(X_{s_2}, \Z(1))
$$ 
mapping $c_1(\l_{s_1})$ to $c_1(\l_{s_2})$, such that $\delta_B(\gamma\cdot x) = \delta_\pi(\gamma)\cdot \delta_B(x)$ for every $x \in H^2_B(X_{s_1}, \Z(1))$ and $\gamma \in \pi_1(V^{\rm an},s_1)$. The isometry $\delta_B$ defines thus an isometry between $P^2_B(X_{s_1}, \Z(1))$ and $P^2_B(X_{s_2}, \Z(1))$ which we will denote again by $\delta_B$. If the level $\Kg$-structure on $\pi \colon X \lr V$ is given by the class $\a$ in $\{\Kg \backslash \textrm{Isometry}\bigl(L_{2d}, P^2(s_1)\bigr)\}^{\pi_1^{\rm alg}(V, s_1)}$ with respect to the geometric point $s_1$ and $\tilde \a$ is a representative of this class, then it is given at the point $s_2$ by the class of $\delta_B \circ \tilde \a$ in $\Kg \backslash \textrm{Isometry}\bigl(L_{2d}, P^2(s_2)\bigr)$. It is $\pi_1^{\rm alg}(V, s_2)$-invariant. See the discussion before Definition 5.1 in Section 5 in \cite{Riz-MK3}. Hence if we take a marking $a \colon H^2_B(X_{s_1},\Z(1)) \lr L_0$ such that $a(c_1(\l_{s_1})) = e_1 - df_1$ then we can take a marking 
$$
 a \circ \delta_B^{-1} \colon H^2_B(X_{s_2},\Z(1)) \lr L_0
$$ 
for which we have that $a(c_1(\l_{s_2})) = e_1 - df_1$. So we see that 
\begin{gather}
 j_{d,\Kg,V}(s_1) = [a\circ h_{s_1} \circ a^{-1}, a \circ \tilde \a]_\Kg\label{connimage1} \\
 j_{d,\Kg,V}(s_2) = [a\circ \delta_B^{-1} \circ h_{s_2} \circ \delta_B \circ a^{-1}, a \circ \delta_B^{-1} \circ \delta_B \circ \tilde \a]_\Kg\label{connimage2}.
\end{gather}
The sheaf $R^2_B\pi_*\Z(1)$ is a local system on $V^{\rm an}$ for every point $s \in V^{\rm an}$ one can find an open neighborhood $V_s$ of $s$ in $V^{\rm an}$ such that the system $R^2_B\pi_*\Z(1)|_{V_s}$ is constant. We can find a marking $a \colon R^2_B\pi_*\Z(1)|_{V_s} \lr (L_0)_{V_s}$ mapping $c_1(\l)$ to $e_1-df_1$. According to Theorem in Section 5 of \cite{Ast-K3} (or Section 9.7 in \cite{Griff-PerIII}) the map
\begin{displaymath}
 j \colon V_s \lr \Omega^\pm
\end{displaymath}
defined by $j(s) = a \circ h_{X^{\rm an}_s} \circ a^{-1}$ is holomorphic. As $V_s$ is connected we may assume that its image in $\Omega^\pm$ under the morphism $j$ is contained in $\Omega^+$. Then we see from \eqref{connimage1} and \eqref{connimage2} that $j_{d,\Kg,V}(V_s) \subset \Gamma_{[g]} \backslash \Omega^+$ where $g = a_s \circ \tilde \a_s$. Further, $pr \colon \Omega^+ \lr \Gamma_{[g]} \backslash \Omega^+$ is holomorphic and we have that $j_{d,\Kg,V}|_{V_s} = pr \circ j$. Hence $j_{d,\Kg,V}|_{V_s}$ is holomorphic.

According to Proposition 3.3.1 in \cite{A-HV} applied to $X^{\rm an}_{V_s} \lr V_{s}$, the holomorphic map $j$ is a local isomorphism and therefore the same holds for $j_{d,\Kg,\C}|_{V_s}$ as $\Omega^+$ is the universal covering space of $\Gamma_{[g]} \backslash \Omega^+$. Those conclusions are valid for a neighborhood of any point $s$ in $V^{\rm an}$ hence we see that $j_{d,\Kg, V} \colon V^{\rm an} \lr \Omega^+\backslash \Gamma$ is holomorphic and it is a local isomorphism.
\newline
\newline
{\bf Step 4:} For every connected component $V$ of $U$ we defined a holomorphic morphism $j_{d,\Kg,V} \colon V^{\rm an} \lr Sh_\Kg(G,\Omega^\pm)_\C$ which is a local isomorphism. By Step 3 it factorizes though a connected component $\Gamma_{[g]} \backslash \Omega^+$ of $Sh_\Kg(G,\Omega^\pm)_\C$ for some $g\in G(\A_f)$ so using a result of A. Borel (Th\'eor\`eme 5.1 in \cite{D-K3}) we conclude that $j_{d,\Kg,V}$ is an algebraic morphism. Indeed, we can apply Th\'eor\`eme 5.1 in {\it loc. cit.} because the group $\Gamma_{[g]}$ is torsion free as $\Kg \subset \Kg_n$ for some $n \geq 3$. We also have that $j_{d,\Kg,V}$, being an analytic local isomorphism, is \'etale. Gluing the morphisms $j_{d,\Kg,V}$ for all connected components $V$ of $U$ we obtain a morphism of $\C$-schemes
$$
 j_{d,\Kg,U} \colon U \lr Sh_\Kg(G,\Omega^\pm)_\C
$$ 
which is \'etale.
\newline
\newline
{\bf Step 5:} We will show that $j_{d,\Kg,U}$ descends to a morphism of algebraic spaces \begin{displaymath}
 j_{d,\Kg,\C} \colon \Fk_{2d,\Kg,\C} \lr Sh_\Kg(G,\Omega^\pm)_\C.
\end{displaymath} 
We have to show that the two projection maps
\begin{displaymath}
 j_{d,\Kg,U} \circ pr_i \colon U \times_{\Fk_{2d,\Kg,\C}} U \lr Sh_\Kg(G,\Omega^\pm)_\C 
\end{displaymath}
for $i = 1,2$ coincide (see Chapter II, Proposition 1.4 in \cite{Knu-AS}). As $\Fk_{2d,\Kg,\C}$ is a reduced algebraic space over $\C$ we have that $U\times_{\Fk_{2d,\Kg,\C}} U$ is a reduced $\C$-scheme (Chapter II, Definition 1.1 in \cite{Knu-AS}). Hence we can check the equality of the two morphisms on $\C$-valued points.

Any $\C$-valued point on $U\times_{\Fk_{2d,\Kg,\C}} U$ is a pair $\big((X_1,\l_1,\a_1),(X_2,\l_2,\a_2),f \bigr)$ where $f$ is an isomorphism of the objects $(X_1,\l_1,\a_1)$ and  $(X_2,\l_2,\a_2)$ in $\Fk_{2d,\Kg,\C}$. Hence from the very definition of the morphism $j_{d,\Kg,U}$ we easily see (just like in the proof of Proposition \ref{Modinterpr1}) that
\begin{displaymath}
j_{d,\Kg,U} \circ pr_1\big((X_1,\l_1,\a_1),(X_2,\l_2,\a_2)\bigr) =
 j_{d,\Kg,U} \circ pr_2 \big((X_1,\l_1,\a_1),(X_2,\l_2,\a_2)\bigr).
\end{displaymath}
Thus we have that $j_{d,\Kg,U} \circ pr_1 = j_{d,\Kg,U} \circ pr_2$ and therefore $j_{d,\Kg,U}$ descends to a morphism $j_{d,\Kg,\C} \colon \Fk_{2d,\Kg,\C} \lr Sh_\Kg(G,\Omega^\pm)_\C$. It is \'etale as $j_{d,\Kg,U}$ is \'etale (Chapter II, Definition 2.1 in \cite{Knu-AS}).
\end{proof}
\begin{cor}
The algebraic space $\Fk_{2d,\Kg,\Q}$ is a scheme.
\end{cor}
\begin{proof}
Combining the above proposition and Corollary 6.17 in Chapter II of \cite{Knu-AS} we conclude that $\Fk_{2d,\Kg,\C}$ is a scheme. Therefore  $\Fk_{2d,\Kg,\Q}$ is a scheme, as well.
\end{proof}
\begin{dfn}\label{PeriodsDef}
The map $j_{d,\Kg,\C}$ is called the \emph{period morphism} (or the \emph{period map}) associated to $d$ and $\Kg$. For every primitively polarized complex K3 surface $(X,\l,\a)$ of degree $2d$ and a level $\Kg$-structure $\a$, the point $j_{d,\Kg,\C}\bigl((X,\l,\a)\bigr) \in Sh_\Kg(G,\Omega^\pm)_\C$ is called the \emph{period point} of $(X,\l,\a)$.
\end{dfn}
\begin{rem}
The period map $j_{d,\Kg,\C}$ defined in the proof of Proposition \ref{PerMorphism} is a slight modification of the period maps used in \cite{Ast-K3} to construct coarse moduli spaces of primitively polarized complex K3 surfaces. We consider moduli spaces over $\Q$ and these have more than one geometric connected component. The morphism constructed above takes this information in to account. We will see later that this is essential for having the period morphism defined over $\Q$. 
\end{rem}
In Section \ref{FinalComm} we show that the image $j_{d,\Kg,\C}(\Fk_{2d,\Kg,\C})$ is dense in $Sh_\Kg(G,\Omega^\pm)_\C$ and that its complement is a divisor.

Suppose that $\Kg$ is an admissible subgroup of $\SO(V_{2d})(\hat \Z)$ (see Definition 5.5 in \cite{Riz-MK3}). Then one can consider the moduli space $\Fk^{\rm full}_{2d,\Kg,\C}$ which is an open subspace of $\Fk_{2d,\Kg,\C}$ (see Theorem 6.8 in \cite{Riz-MK3}). Hence we have a period map
$$
 j_{d,\Kg,\C} \colon \Fk_{2d,\Kg,\C}^{\rm full} \lr  Sh_\Kg(G,\Omega^\pm)_\C
$$
given by the restriction of $j_{d,\Kg,\C} \colon \Fk_{2d,\Kg,\C} \lr Sh_\Kg(G,\Omega^\pm)_\C$ to $\Fk_{2d,\Kg,\C}^{\rm full}$. We will show below that this period morphism is injective. This result is a direct consequence of the \emph{global Torelli theorem for K3 surfaces} of Piatetskij-Shapiro and Shafarevich. See Expos\'e VII, Section 3 in \cite{Ast-K3} and Corollary 11.2 in~\cite{BPvV}.
\begin{prp}\label{InjPerMorphism}
For an admissible subgroup $\Kg$ of $\SO(V_{2d})(\hat \Z)$ the period map 
$$
 j_{d,\Kg,\C} \colon \Fk_{2d,\Kg,\C}^{\rm full} \lr  Sh_\Kg(G,\Omega^\pm)_\C
$$ 
is an open immersion.
\end{prp} 
\begin{proof}
Note first that the algebraic space $\Fk_{2d,\Kg,\C}^{\rm full}$ is a scheme as we have an open immersion $i_\Kg \colon \Fk_{2d,\Kg}^{\rm full} \hookrightarrow \Fk_{2d,\Kg}$ into a scheme (cf. Theorem 6.8 in \cite{Riz-MK3}). Further, by Proposition \ref{PerMorphism} the map $j_{d,\Kg,\C} \colon \Fk_{2d,\Kg,\C}^{\rm full} \lr  Sh_\Kg(G,\Omega^\pm)_\C$ is \'etale hence it is open. We have to show that it is an immersion. As both schemes are reduced it is enough to show that the morphism is injective on $\C$-valued points.

Suppose that $(X_i,\l_i,\a_i) \in \Fk_{2d,\Kg,\C}^{\rm full}(\C)$ for $i = 1,2$ are two points such that 
$$
 j_{d,\Kg,\C}\bigl((X_1,\l_1,\a_1)\bigr) =  j_{d,\Kg,\C}\bigl((X_2,\l_2,\a_2)\bigr)
$$ 
in $Sh_\Kg(G,\Omega^\pm)_\C$. Let $\tilde \a_i$ be two representatives of the classes $\a_i$ and let 
$$
a_i \colon H^2_B(X_i,\Z(1)) \lr L_0
$$ 
be two isometries as in the definition of the map $j_{d,\Kg,\C}$. Then we have that 
$$
[a_1\circ h_{X_1} \circ a_1^{-1}, a_1 \circ \tilde \a_1]_\Kg = [a_2\circ h_{X_2} \circ a_2^{-1}, a_2 \circ \tilde \a_2]_\Kg.
$$ 
Hence there are two elements $q \in G(\Q)$ and $\kappa \in \Kg$ such that we have an equality
\begin{equation}\label{InjectjFull}
\bigl(q\cdot (a_1\circ h_{X_1} \circ a_1^{-1}), q\circ a_1 \circ \tilde \a_1 \circ \kappa\bigr) =
\bigl(a_2\circ h_{X_1} \circ a_2^{-1}, a_2 \circ \tilde \a_2\bigr).
\end{equation}
From the equality between the second elements in \eqref{InjectjFull} we see that $q = a_2\circ \tilde \a_2 \circ \kappa \circ \tilde \a_1^{-1} \circ a_1^{-1}$. Hence it belongs to $\{g \in \SO(V_{2d})(\hat \Z)\ |\ g(e_1-df_1) = e_1-df_1\}$, by the very definition of a full level $\Kg$-structure, and being in $G(\Q)$ we conclude that $q \in \{g \in \SO(V_{2d})(\Z)\ |\ g(e_1-df_1) = e_1-df_1\}$. The equality between the first elements in \eqref{InjectjFull} shows that 
$$
 a_2^{-1} \circ q \circ a_1 \colon H^2_B(X_2,\Z(1)) \lr H^2_B(X_1,\Z(1))
$$ 
is a Hodge isometry, mapping the class of $\l_2$ to the class of $\l_1$ and preserving the level structures. By the global Torelli theorem for K3 surfaces one concludes that it comes from an isomorphism of the triples $(X_1,\l_1,\a_1)$ and $(X_2,\l_2,\a_2)$. Therefore the morphism $j_{d,\Kg, \C}$ is an immersion.
\end{proof}
\begin{rem}\label{NonInjPerMap}
In general, the morphism $j_{d,\Kg,\C} \colon \Fk_{2d,\Kg,\C} \lr Sh_\Kg(G,\Omega^\pm)_\C$ need not be injective. We cannot apply the arguments of the proof of Proposition \ref{InjPerMorphism} as we only get a Hodge isometry
$$
 a_2^{-1} \circ q \circ a_1 \colon P^2_B(X_2,\Z(1)) \lr P^2_B(X_1,\Z(1)).
$$ 
As not every such isometry is induced by an isometry between the cohomology groups $H^2_B(X_2,\Z(1))$ and $H^2_B(X_1,\Z(1))$, mapping $c_1(\l_2)$ to $c_1(\l_1)$, we cannot conclude that $(X_1,\l_1)$ and $(X_2,\l_2)$ are isomorphic.
\end{rem}
\section{Complex Multiplication for K3 Surfaces}
Here we will prove that the field of definition of $j_{d,\Kg,\C}$ is $\Q$. This is an analogue of Th\'eor\`eme 4.21 in \cite{D-ShV1} concerning periods of abelian varieties. We will do this first by proving a variant of the main theorem of complex multiplication for abelian varieties (Theorem 5.3 in Chapter I of \cite{Mil-CanModMixShV}) in the case of exceptional K3 surfaces and then applying a density result for those surfaces. We will carry out this strategy in Sections \ref{CMK3Stat}-\ref{CMK3SSect}. Before that we give a short review of some results from class field theory and canonical models of Shimura varieties.

We begin by making the following notation which will be used from now on.

Let $X/\C$ be a non-singular projective variety and consider an automorphism $\sigma \in \Aut(\C)$. Let $X^\sigma$ be conjugate of $X$ by $\sigma$. For any $n \in \Z$ we will denote by 
\begin{displaymath}
\sigma_{X,f} \colon H^i_{\rm et}(X,\A_f(n)) \lr H^i_{\rm et}(X^\sigma, \A_f(n))
\end{displaymath}
or simply by $\sigma_f$, the morphism on \'etale cohomology induced by $\beta$.

For a non-singular projective surface $X$ the morphism $\beta$ induces a morphism 
$$
 \beta^* \colon \Pic(X) \lr \Pic(X^\sigma)
$$
and we will denote it by $\sigma_{\Pic}$. Recall that we have a decomposition 
$$
 H^2_{\rm et}(X,\A_f(1))= A_{X,\A_f} \oplus T_{X,\A_f}
$$
and similarly
$$
 H^2_{et}(X^\sigma,\A_f(1))= A_{X^\sigma,\A_f} \oplus T_{X^\sigma,\A_f}.
$$
We have that
$$
\sigma_{X,f} = \sigma_{\Pic, \A_f} \oplus \sigma_{X,f}|_{T_{X,\A_f}}
$$
where $\sigma_{\Pic, \A_f}$ is the morphism sending $c_1(\l)$ to $c_1(\l^\sigma)$ for any $\l \in \Pic(X)$. In the sequel, we shall use the notation $\sigma_{X,f}$ for the morphism $\sigma_{X,f}|_{T_{X\A_f}} \colon T_{X,\A_f} \lr T_{X^\sigma, f}$.
\subsection{Class Field Theory}\label{ClassFieldTheory}
Let $E$ be a number field and denote by $E^{\rm ab}$ the maximal abelian extension of $E$. Class field theory provides us with a description of $\Gal(E^{\rm ab}/E)$. There exists a surjective homomorphism
\begin{displaymath}
\text{rec}_E \colon \A_E^\times \lr \Gal(E^{\rm ab}/E)
\end{displaymath}
such that $E^\times$ is in its kernel and for every finite abelian extension $L$ of $E$ the following diagram
$$
\xymatrix{
E^\times\backslash \A^\times_E \ar[d] \ar[rr]^{\text{rec}_E}_{\text{onto}} & &
\Gal(E^{\rm ab}/E) \ar[d]^{\sigma \mapsto \sigma|_L} \\
E^\times \backslash \A_E^\times / \Nm_{L/E}(\A_L^\times)
\ar[rr]^{\text{rec}_{L/E}}_{\cong} & & \Gal(L/E)
}
$$
is commutative. We refer to Chapter 3 in \cite{N-CFT} and Chapter VII in \cite{CF-CFT} for proofs and some properties of this homomorphism. To make notations easier when considering canonical models of Shimura varieties we define the map
\begin{displaymath}
 \text{art}_E \colon \A^\times_E \lr \Gal(E^{\rm ab}/E)
\end{displaymath}
by $\text{art}_E(\a) = \text{rec}_E(\a)^{-1}$.
\subsection{The Homomorphism $r_h$}
Let $V$ be a finite dimensional $\Q$-vector space and let
$$
 h \colon \s \lr \GL(V)_\R
$$
be a $\Q$-HS on $V$. Let $T \subset \GL(V)$ be a $\Q$-torus and suppose that the homomorphism $h$ factorizes thought $T_\R$. Then the same holds for the cocharacter $\mu_h$ (cf. Section \ref{notatch3}) and we have that
$$
 \mu_h \colon \Gm_{m,\C} \lr T_\C
$$
is defined over $\bar \Q$. Let $E(h)$ be the field of definition of $\mu_h$ i.e., the reflex field of the pair $(T,h)$. It is a number field. Composing $\mu_h \colon \Gm_{m,E(h)} \lr T_{E(h)}$ with he norm morphism we obtain a homomorphism
\begin{displaymath}
 r(T,h) \colon \text{Res}_{E(h)/\Q}(\Gm_{m,E(h)}) \lr T.
\end{displaymath}
\begin{dfn}\label{rhm}
With notations as above define the homomorphism $r_h \colon \A_{E(h)}^\times \lr T(\A_f)$ as being the composition
$$
\xymatrix{
r_h \colon \A_{E(h)}^\times \ar@{=}[r] & \text{Res}_{E(h)/\Q}(\Gm_{m,E(h)})(\A) \ar[rrr]^{r(T,h)} & & & T(\A) \ar[r]^{\text{proj}} & T(\A_f).
}
$$
\end{dfn}
\subsection{The Canonical Model of $Sh_\Kg(G,\Omega^\pm)_\C$}\label{CanModSect}
Let $\Kg$ be a compact open subgroup of $G(\A_f)$. The canonical model $Sh_\Kg(G,\Omega^\pm)$ of $Sh_\Kg(G,\Omega^\pm)_\C$ is scheme over $\Q$ (which is the reflex field of the Shimura datum $(G,\Omega^\pm)$) such that:
\begin{enumerate}
\item [(i)]  one has an isomorphism $Sh_\Kg(G,\Omega^\pm) \otimes_\Q \C \lr Sh_\Kg(G,\Omega^\pm)_\C$;
\item [(ii)] $\Aut(\C)$ acts on $Sh_\Kg(G,\Omega^\pm)_\C$ via the isomorphism given by (i) as follows: For every special pair $(T,x)$ of $(G,\Omega^\pm)$ one has that
$$
 \sigma[x,a]_\Kg = [x,r_x(s)a]_\Kg
$$
for all $\sigma \in \Aut(\C/E(x))$ and $s \in \A_{E(x)}^\times$ such that ${\rm art}_{E(x)}(s) = \sigma|_{E(x)^{\rm ab}}$. Here the morphism $r_x \colon \A_{E(x)}^\times \lr T(\A_f)$ is the one associated to the pair $(T,x)$ as in Definition \ref{rhm}.
\end{enumerate}
These two properties determine the scheme $Sh_\Kg(G,\Omega^\pm)$ uniquely up to a unique isomorphism. For details concerning canonical models of Shimura varieties and their properties we refer to Section 2 in \cite{BM-SV}.
\subsection{Statement of the Main Theorem of Complex Multiplication for Exceptional K3 Surfaces}\label{CMK3Stat}
Let $X$ be an exceptional K3 surface of CM-type $(E_X,\e_X)$ over $\C$. As in the case of abelian varieties with complex multiplication we are interested in a relation between the various cohomology groups of $X$ and its conjugate $X^\sigma$ by an automorphism $\sigma$ of $\C$. In this section we will state the main results of complex multiplication for exceptional K3 surfaces. To make notations easier we will denote by 
$$
 E:= \e_X(E_X) \subset \C
$$
the reflex field of $\Hg(X)$.

Recall that the Hodge structure homomorphism $h_X \colon \s \lr \SO(P^2_B(X,\R(1)))$ factorizes
$$
 h_X \colon \s \lr \Hg(X)_\R \subset \SO(T_{X,\R}) \hookrightarrow \SO(P^2_B(X,\R(1))).
$$
Let $\mu_X \colon \Gm_{m,E} \lr \Hg(X)_E$ be the corresponding cocharacter and let 
$$
 r_X \colon \A^\times_E \lr \Hg(X)(\A_f) \subset \SO(T_{X,\Q})(\A_f)
$$ 
be the morphism associated to $(\Hg(X),h_X)$ as in Definition \ref{rhm}.
\begin{lem}\label{RefFTwist}
Suppose given an exceptional K3 surface $X$ of CM-type $(E_X,\e_X)$. If $\sigma \in \Aut(\C/E)$, then $X^\sigma$ is an exceptional K3 surface and the reflex field of $\Hg(X^\sigma)$ is~$E$.
\end{lem}
One can give a proof of the lemma using a `a Hodge cycle is an absolute Hodge cycle' argument. We will give a proof in Section \ref{TheProof} using abelian surfaces.

Let $X$ be an exceptional K3 surface of CM-type $(E_X,\e_X)$ and let $\sigma \in \Aut(\C/E)$. Then by a \emph{$E_X$-linear isometry} $\eta \colon T_{X,\Q} \lr T_{X^\sigma,\Q}$ we shall mean an isometry $\eta$ such that
$$
 \eta(e\cdot t) = (\e_{X^\sigma}^{-1} \circ \e_X)(e) \cdot \eta(t)
$$
for every $t \in T_{X,\Q}$ and $e \in E_X$.
\begin{thm}[Complex multiplication for exceptional K3 surfaces]\label{MainThCMEK3}
Let $X$ be an exceptional K3 surface of CM-type $(E_X,\e_X)$. Let $E = e_X(E_X) \subset \C$ be its reflex field and let $\sigma \in \Aut(\C/E)$. Then for any id\`ele $s\in \A^\times_{E}$ with ${\rm art}_{E}(s) =
\sigma|_{E^{\rm ab}}$ there is a unique $E_X$-linear isomorphism of polarized $\Q$-HS
$$
 \eta_X \colon T_{X,\Q} \lr T_{X^\sigma,\Q}
$$
such that $\eta_{X,f}(r_X(s)t) = \sigma_{X,f}(t)$ for every $t \in T_{X,\A_f}$.
\end{thm}
If such $\eta_X$ exists, then it is necessarily unique. Indeed, the condition imposed on $\eta_{X,f}$ determines it uniquely and hence $\eta$ is also determined uniquely via the natural isomorphism $T_{X,\A_f} \cong T_{X,\Q} \otimes \A_f$.

We will give a proof of the theorem in Section \ref{TheProof}. We will use first a geometric construction due to Shioda and Inose to reduce the problem to a similar statement for abelian surfaces. Then we will show that the corresponding statement for abelian surfaces follows from the main theorem of complex multiplication for abelian varieties. We present these results in the next two sections.
\begin{rem}
We wonder if one could give a `direct' proof of Theorem \ref{MainThCMEK3} similar to the proof of Theorem 11.2 in \cite{Mil-IShV} using, for instance, arguments of the type `a Hodge cycle is an absolute Hodge cycle' on a K3 surface. This can be done in the case $(X,\l)$ is defined over an intermediate field $E \subset K \subset \C$ and $\sigma \in \Aut(\C/K)$.
\end{rem}
\subsection{The Results of Shioda and Inose}\label{ShiodaInoseResults}
We shall describe a geometrical way for constructing exceptional K3 surfaces with given transcendental lattice using product abelian surfaces. We will follow the exposition of Shioda and Inose in their paper \cite{Sh-I} with some notational differences.

Let $A = C_1 \times C_2$ be a product of two elliptic curves over $\C$ and let $Y$ be the Kummer surface associated to $A$. Let $\pi \colon \tilde A \lr A$ be the the blowing up of the 2-torsion of $A$ and let $[-1]_{\tilde A}$ be the involution on $\tilde A$ induced by the automorphism $[-1]_A$ of $A$. Denote by $\iota \colon \tilde A \lr \tilde A/\langle{[-1]_{\tilde A}}\rangle = Y$ be the quotient morphism of degree 2. One has morphisms induced on Betti cohomology with $\Z$-coefficients and hence on the corresponding transcendental lattices
$$
 \pi^* \colon T_A \lr T_{\tilde A} \ \ \ \text{and} \ \ \ \iota^* \colon T_Y \lr T_{\tilde A}.
$$
We know that $\pi^*$ is an isomorphism of polarized $\Z$-HS and $\iota^*$ is an isomorphism of $\Z$-HS multiplying the intersection form by 2 i.e., $(\iota^*x,\iota^*y)_{\tilde A} = 2(x,y)_Y$.

Let $\{u_i\}_{i=1}^4$, $\{v_j\}_{j=1}^4$ be the four $2$-torsion points on $C_1$ and $C_2$ respectively. Denote by $E_{ij}$ the sixteen non-singular rational curves on $Y$
corresponding to the points $(u_i,v_j)$ on $A$. In other words we have that $E_{ij} = \iota\bigl(\pi^{-1}(u_i,v_j)\bigr)$. Following the notations of \cite{Sh-I} we denote by $F_i$ and $G_j$ the non-singular rational curves $\iota\bigl(\pi^{-1}(u_i\times C_2)\bigr)$ and $\iota\bigl(\pi^{-1}(C_1\times v_j)\bigr)$ on $Y$.

Consider the divisor
$$
 D = E_{21} + 2F_2 + 3E_{23} +4G_3 + 5E_{13} +6F1 + 3E_{12} + 4E_{14} + 2G_4
$$
on $Y$. By Lemma 1.1 in \cite{Sh-I} the linear system $|D|$ gives a morphism $\Phi \colon Y \lr \P^1$ of which $D$ is a singular fiber, say $D = \Phi^{-1}(t_0)$ for some $t_0 \in \P^1$. We look further at two divisors
$$
 B_1 = F_3 + E_{31} + E_{32} \ \ \ \text{and} \ \ \  B_2 = F_4 + E_{41} + E_{42}
$$
on $Y$. One can see that they do not meet $D$ and their supports are connected. Hence we conclude that the image $\Phi(B_i)$ is a point $t_i$ in $\P^1$ and $B_i$ is contained in the singular fiber $\Phi^{-1}(t_i)$ for $i=1,2$ (see the figure on page 122 of \cite{Sh-I} and the comments following it). Let $f \colon \P^1 \lr \P^1$ be the finite morphism of degree 2 branched only at $t_1$ and $t_2$ and consider the fiber product $Y\times_{\Phi, \P^1, f} \P^1$.
\begin{lem}
The surface $Y \times_{\P^1} \P^1$ has a minimal model $X$ which is a K3 surface (hence it is unique).
\end{lem}
\begin{proof}
See Lemma 3.1 in \cite{Sh-I}.
\end{proof}
The elliptic pencil $\Phi \colon Y \lr \P^1$ on $Y$ induces an elliptic pencil $\Psi \colon X \lr \P^1$ on $X$ (see \S 3, p. 124 in \cite{Sh-I}). The morphism $f \colon \P^1 \lr \P^1$ induces an involution of the surface $Y\times_{\Phi, \P^1, f} \P^1$. Therefore it induces an involutive birational transformation of $X$, hence by the minimality of a K3 surface an automorphism $a$ of $X$. It has 8 fixed points $\{p_i\}_{i=1}^8$ and $Y$ is the minimal model of the quotient surface $X/a$. For details see \S 8, pages 585-586, 591-592, 600-602 in \cite{Kod-CASII} and the remarks after the proof of Lemma 3.1 on page 125 in~\cite{Sh-I}.

Let $\beta \colon \tilde X \lr X$ be the blow-up of the 8 points $p_i$, $i=1,\dots,8$, on $X$. Then the involution $a$ on $X$ induces an involution $\tilde a$ on $\tilde X$. If we denote the quotient morphism $\tilde X \lr \tilde X/\tilde a$ by $\gamma$ then we have the following commutative diagram
$$
\xymatrix{
& {\tilde X} \ar[dl]_\gamma \ar[dr]^\beta & \\
Y = \tilde X/{\tilde a} \ar[dr] & & X. \ar[dl] \\
& X/a &
}
$$
The degree of the morphism $\gamma$ is 2. The map $\beta$ induces an isomorphism of polarized $\Z$-HS $\beta^* \colon T_X \lr T_{\tilde X}$. The main result of Section 2 of \cite{Sh-I} is that $\gamma^* \colon T_Y \lr T_{\tilde X}$ is an isomorphism of $\Z$-HS such that $(\gamma^*x,\gamma^*y)_{\tilde X} = 2(x,y)_Y$. Putting the preceding two diagrams together we obtain
\begin{equation}\label{ShIConstruction}
\xymatrix{
& \tilde A \ar[dl]_\pi \ar[dr]^\iota & & \tilde X \ar[dl]_\gamma \ar[dr]^\beta & \\
A & & Y & & X.
}
\end{equation}
Shioda and Inose describe the relation between the transcendental lattices of $A$ and $X$ using those morphisms.
\begin{thm}[Shioda-Inose]\label{ShIEK3-1}
With the notations as above one has that the morphism
$$
 \phi \colon T_X \lr T_A.
$$
defined as $\phi = \pi^{* -1} \circ \iota^* \circ \gamma^{* -1} \circ \beta^*$ induces an isomorphism of polarized $\Z$-HS.
\end{thm}
\begin{proof}
The main difficulty is to prove that the map $\gamma^*$ is an isomorphism. We refer to the proof of Theorem 2 in \cite{Sh-I}. Note that Shioda and Inose use homology groups and we use cohomology groups. But in our case all those groups are free and we obtain the result using duality.
\end{proof}
\begin{rem}
Note that a priori the whole construction depends on choosing a numbering of $A[2](\C)$. We shall be interested in constructing exceptional K3 surfaces. As we will see below for these surfaces the choices involved change only the morphisms $\beta$ and $\gamma$ but not the surface $X$ itself.
\end{rem}
\begin{rem}
Note that by the comparison theorem between Betti and \'etale cohomology the map $\phi_f = \pi^{* -1}_f \circ \iota^*_f \circ \gamma^{* -1}_f \circ \beta^*_f$ induces an isomorphism
$$
 \phi_f \colon T_{X,\A_f} \lr T_{A,\A_f}.
$$
Indeed, we have that $\phi_f = \phi \otimes_\Z \A_f$ and we know that $\phi$ is an isomorphism.
\end{rem}
In order to explain the construction in the proof of the main result of \cite{Sh-I} we will follow their notations working with homology instead of cohomology. If $X$ is a non-singular projective surface over $\C$ we will denote by $T^{\rm hom}_X$ the homological transcendental lattice. In other words we define $T^{\rm hom}_X = (\Pic(X))^\perp \subset H_2(X,\Z(-1))$.

Let $X$ be an exceptional K3 surface over $\C$. Denote by $p_X$ the period on
$T^{\rm hom}_X$ i.e., the linear functional, determined up to a constant by
$$
 p_X(t) = \int_t\omega_X
$$
for $t\in T_X^{\rm hom}$ and $\omega_X$ a non-vanishing holomorphic 2-form on $X$. We say
that a basis $\{y_1,y_2\}$ is oriented if the imaginary part of
$p_X(y_1)/p_X(y_2)$ is positive.

Let $\{y_1,y_2\}$ be an oriented basis of $T_X^{\rm hom}$. In it the bilinear form on $T_X^{\rm hom}$ is given by a matrix
\begin{equation}\label{TrFormMatrix}
 Q =\biggl(
\begin{matrix}
\langle y_1, y_1\rangle & \langle y_1,y_2 \rangle \\
\langle y_2,y_1\rangle & \langle y_2, y_2 \rangle
\end{matrix}\biggr)
=
\biggl(
\begin{matrix}
2a & b \\
b & 2c
\end{matrix}
\biggr)
\end{equation}
for some $a,b,c \in \Z$ with $a,c > 0$ and $\Delta = b^2 - 4ac < 0$. Then we have that $E_X = \End_{\rm HS}(T_{X,\Q}) = \End_{\rm HS}(T_{X,\Q}^{\rm hom})$ is isomorphic to $\Q(\sqrt{\Delta}) \subset \C$. In our notations from Section \ref{CMK3Stat} we have that $e_X(E) = E = \Q(\sqrt{\Delta})$.

Let $\tau_1 = (-b + \sqrt{\Delta})/2a$ and $\tau_2 = (b +\sqrt{\Delta})/2$ and consider the elliptic curve $C_i = \C/\Lambda_{\tau_i}$ where $\Lambda_{\tau_i} = \Z + \Z\tau_i$. These elliptic curves are isogenous and have complex multiplication by $E$.
\begin{thm}[Shioda-Inose]\label{ShIEk32}
Let $X$ be an exceptional K3 surface over $\C$ and consider the product abelian surface $A = C_1\times C_2$ where $C_i$ for $i=1,2$ are the CM elliptic curves defined above. Then the K3 surface $X_A$ constructed in Theorem \ref{ShIEK3-1} is isomorphic to $X$.
\end{thm}
\begin{proof}
We refer to the proof of Theorem 4 in \cite{Sh-I}. The idea is to compare the lattices $T^{\rm hom}_X$ and $T^{\rm hom}_{X_A}$. Using Theorem \ref{ShIEK3-1} one sees that those two lattices are isometric. A result of Piatetskij-Shapiro and Shafarevich says that an exceptional K3 surface is uniquely determined by its transcendental lattice (see \S 8 in \cite{PSh-Sh} and also the remarks made in \cite{ShM-TorAS}). Hence one concludes that $X$ is isomorphic to $X_A$.
\end{proof}
\begin{rem}
Note that if $X$ is an exceptional K3 surface then the construction described in Theorem \ref{ShIEk32} is independent of the numbering of $A[2](\C)$. Indeed, starting with a numbering of $A[2](\C)$ one can constructs an exceptional K3 surface $X_1$ and an isomorphism of polarized $\Z$-HS $\phi_1 \colon T_{X_1} \lr T_A$. Starting with a different numbering and different $f$ one constructs an exceptional K3 surface $X_2$ with an isomorphism of polarized $\Z$-HS $\phi_2 \colon T_{X_2} \lr T_A$. Hence $T_{X_1}$ and $T_{X_2}$ are isometric and by the result of Piatetskij-Shapiro and Shafarevich $X_1$ and $X_2$ are isomorphic. Note that the morphisms involved in the construction might change.
\end{rem}
\subsection{Complex Multiplication for Product Abelian Surfaces}
Let $A$ be a complex abelian surface of CM type $(E,\Phi)$ and denote by $E^*$ its reflex field. Let $\sigma$ be an element of $\Aut(\C/E^*)$. The main theorem of complex multiplication gives a relation between Betti and \'etale cohomology of $A$ and $A^\sigma$. We will need this in a special case. Before stating the result we introduce some notations.

Let $E \subset \C$ be a quadratic imaginary field and let $C_1$ and $C_2$ be two elliptic curves with CM by $E$. One should keep in mind here the data of Theorem \ref{ShIEk32}. Then $C_i$ is of CM-type $E \subset \C$. Let $A$ be the product abelian surface $C_1\times C_2$. Then the reflex field of the torus $\MT(A)$ is $E \subset \C$ (see Chapter IV, \S 18.7 in \cite{Sh-CMAV}). Consider the transcendental space $T_{A,\Q}$ and define
$$
 E_A := \End_{\rm HS}(T_{A,\Q}).
$$ 
Then $E_A$ is a quadratic imaginary field. The reflex field of the torus $\MT(T_{A,\Q})$ is $E \subset \C$. On the other hand if $\e_A \colon E_A \lr \End_\C(H^{2,0}(A)) \cong \C$ denote the action of $E_A$ on the space of holomorphic two-forms on $A$, then just like in the case of K3 surfaces the field $\e_A(E_A) \subset \C$ is the reflex field of $\MT(T_{A,\Q})$. Hence we have an isomorphism $\e_A \colon E_A \lr E$.

We have natural isomorphisms of cohomology groups
\begin{equation}\label{cohisom1}
H^1_B(A,\Z) \cong H^1_B(C_1,\Z) \oplus H^1_B(C_2,\Z) \ \ \  \text{and} \ \ \ H^1_{\rm et}(A,\hat \Z) \cong H^1_{\rm et}(C_1,\hat \Z) \oplus H^1_{\rm et}(C_2,\hat \Z).
\end{equation}
If $h \colon \s \lr \GL(H^1_B(A,\R))$ and $h_i \colon \s \lr \GL(H^1_B(C_i,\R))$, for $i=1,2$ are the corresponding homomorphisms defining the three $\Z$-HS, then we have that $h = h_1 \oplus h_2$. Hence we have that $\mu_h = \mu_{h_1} \oplus \mu_{h_2}$ where $\mu_h$ and $\mu_{h_i}$ are the cocharacters defined in Section \ref{notatch3}. Further we know that
\begin{equation}\label{cohisom2}
 H^2_B(A, \Z(1)) \cong \bigl(\wedge^2 H^1_B(A,\Z)\bigr) \otimes \Z(1) \ \ \ \text{and} \ \ \ H^2_{\rm et}(A, \hat \Z(1)) \cong \bigl(\wedge^2 H^1_{\rm et}(A,\hat \Z)\bigr) \otimes \hat \Z(1)
\end{equation}
and therefore combining \eqref{cohisom1} and \eqref{cohisom2} we have natural isomorphisms
\begin{equation}\label{cohisom3}
\begin{matrix}
H^2_B(A,\Z(1)) \cong \\
\biggl(\wedge^2 H^1_B(C_1,\Z)\otimes \Z(1)\biggr) \oplus \biggl(\wedge^2 H^1_B(C_2,\Z) \otimes \Z(1)\biggr) \oplus \biggl(H^1_B(C_1,\Z) \otimes H^1_B(C_2,\Z) \otimes \Z(1)\biggr)
\end{matrix}
\end{equation}
and
\begin{equation}\label{cohisom4}
\begin{matrix}
H^2_{\rm et}(A,\hat \Z(1)) \cong \\
\biggl(\wedge^2 H^1_{\rm et}(C_1,\hat \Z)\otimes \hat \Z(1)\biggr) \oplus \biggl(\wedge^2 H^1_{\rm et}(C_2,\hat \Z) \otimes \hat \Z(1)\biggr) \oplus \biggl(H^1_{\rm et}(C_1,\hat \Z) \otimes H^1_{\rm et}(C_2,\hat \Z) \otimes \hat \Z(1)\biggr).
\end{matrix} 
\end{equation}
The spaces $\bigl(\wedge^2 H^1_B(C_1,\Q)\bigr)\otimes \Q(1)$ and $\bigl(\wedge^2 H^1_B(C_1,\Q)\bigr)\otimes \Q(1)$ (respectively with $\A_f$-coefficients) consist of algebraic classes. Hence for the homomorphism 
$$
 h_A \colon \s \lr \GL(H^2_B(A,\R(1)))
$$ 
giving the $\Z$-HS on $H^2_B(A,\Z(1))$ we have 
\begin{equation}
h_A = (\wedge^2 h) \otimes h_{\Z(1)} = \bigl(\wedge^2 h_1 \otimes h_{\Z(1)}\bigr) \oplus \bigl(\wedge^2 h_2 \otimes h_{\Z(1)}\bigr) \oplus \bigl(h_1\otimes h_2 \otimes h_{\Z(1)}\bigr).
\end{equation} 
Then for the corresponding cocharacters one has 
\begin{equation}
\mu_A = (\wedge^2 \mu_h) \otimes \mu_{\Z(1)} = \bigl(\wedge^2 \mu_1 \otimes \mu_{\Z(1)}\bigr) \oplus \bigl(\wedge^2 mu_2 \otimes \mu_{\Z(1)}\bigr) \oplus \bigl(\mu_1\otimes \mu_2 \otimes \mu_{\Z(1)}\bigr).
\end{equation}
As we explained in Section \ref{notatch3} the homomorphism $h_A$ and the cocharacter $\mu_A$ factor trough $\SO(T_{A,\Q})$. The Mumford-Tate group $\MT\bigl(H^2_B(A,\Q(1))\bigr)$ is a torus and we have homomorphisms of algebraic groups
$$
 h_A \colon \s \lr \MT\bigl(H^2_B(A,\Q(1))\bigr)_\R \subset \SO(T_{A,\R})
$$
and
$$
 \mu_A \colon \Gm_{m,\C} \lr \MT\bigl(H^2_B(A,\Q(1))\bigr)_\C \subset \SO(T_{A,\C}).
$$
We have that $\MT(H^2_B(A,\Q(1))) = \MT(T_{X,\Q})$. The field of definition of $\mu_A$ is $E \subset \C$. Let
$$
 r_A \colon \A_{E,f}^\times \lr \MT(T_{X,\Q})(\A_f) \subset \SO(T_{A,\Q})(\A_f)
$$
be the morphism associated to $(\MT(T_{X,\Q}),h_A)$ as in Definition \ref{rhm}.

Let $\sigma \in \Aut(\C/E)$. Then by a \emph{$E_A$-linear isometry} $\eta \colon T_{A,\Q} \lr T_{A^\sigma,\Q}$ we shall mean an isometry $\eta$ such that
$$
 \eta(e\cdot t) = (\e_{A^\sigma}^{-1} \circ \e_A)(e) \cdot \eta(t)
$$
for every $t \in T_{A,\Q}$ and $e \in E_X$. Note that this definition is correct as the reflex fields of $\MT(A)$ and $\MT(A^\sigma)$ are $E \subset \C$.
\begin{prp}\label{CMPrAS}
Let $A = C_1\times C_2$ be a product of two elliptic curves with CM by a
quadratic imaginary field $E$. Let $\sigma$ be in $\Aut(\C/E)$ and let $s\in
\A_{E}^\times$ be an id\`ele such that ${\rm art}_E(s) = \sigma|_{E^{\rm ab}}$.
Then there exists an isogeny $\eta \colon A^\sigma \lr A$ such that for the isometry $\eta_f^* \colon T_{A,\A_f} \lr T_{A^\sigma,\A_f}$ induced by $\eta$ acting on \'etale cohomology we have that $\eta_f^*(r_A(s)t) = \sigma_f(t)$ for every $t \in T_{A,\A_f}$.
\end{prp}
\begin{proof}
By Theorem 11.2 in \cite{Mil-IShV} we can find two isogenies $\eta_i \colon C_i^\sigma \lr C_i$ for $i=1,2$ such that for the maps
$$
\eta_{i,f}^* \colon H^1_{\rm et}(C_i,\Q) \lr H^1_{\rm et}(C_i^\sigma,\Q)
$$
we have that $\eta^*_{i,f}(r_i(s)t) = \sigma_{C_i,f}(t)$ for every $t \in H^1_{\rm et}(C_i,\A_f)$, for $i=1,2$. Here 
$$
 r_i \colon \A_{E}^\times \lr \MT(C_i)(\A_f) \hookrightarrow \GL(H^1_B(C_i,\Q))(\A_f)
$$ 
is the homomorphism associated to $(\MT(C_i),h_i)$.

Let $\eta \colon A^\sigma \lr A$ be the product isogeny $(\eta_1,\eta_2)$. It defines an  $E_A$-linear isometry $\eta^*_\Q \colon T_{A,\Q} \lr T_{A^\sigma,\Q}$. Using the decompositions 
\begin{displaymath}
 H^2_{\rm et}(A,\A_f(1)) = T_{A,\A_f} \oplus A_{A,\A_f}\ \ \  \text{and}\ \ \ H^2_{\rm et}(A^\sigma,\A_f(1)) = T_{A^\sigma,\A_f} \oplus A_{A^\sigma,\A_f}
\end{displaymath}
we see that $\eta_f^*|_{A_{A,\A_f}} \colon A_{A,\A_f} \lr A_{A^\sigma,\A_f}$ sends a class $c_1(\l)$ for $\l \in \Pic(A)$ to $c_1(\l^\sigma)$. Further, using the natural isomorphisms \eqref{cohisom1}, \eqref{cohisom2}, \eqref{cohisom3} and \eqref{cohisom4} we see that for 
$$
 \eta^*_f \colon T_{A,\A_f} \lr T_{A^\sigma,\A_f}
$$ 
we have that $\eta_f^*(r(s)t) = \sigma_f(t)$. Here $r\colon \A_E^\times \lr \MT(T_{X,\Q})(\A_f)$ is the morphism obtained as in Definition \ref{rhm} using the cocharacter $\bigl(\wedge^2 (\mu_1\oplus \mu_2)\bigr) \otimes \mu_{\Q(1)}$ which is exactly $\mu_A$. So we have that $r = r_A$ which finishes the proof.
\end{proof}
\subsection{Proof of the Main Theorem of Complex Multiplication for Exceptional K3 Surfaces}\label{TheProof}
Here we will give a proof of the claims announced in Section \ref{CMK3Stat}. We begin by putting together the results of the previous two sections.

Let $X$ be an exceptional K3 surface over $\C$ of CM-type $(E_X,\e_X)$. Let 
$$
 E = \e_X(E_X) = \Q(\sqrt{\Delta})
$$ 
be the quadratic imaginary field defined by the discriminant of the form \eqref{TrFormMatrix} in Section \ref{ShiodaInoseResults}. Let $A$ be the product abelian surface as in Theorem \ref{ShIEk32} associated to $X$. Let us further set $E_A = \End_{\rm HS}(T_{A,\Q})$. The two fields $E_X$ and $E_A$ are isomorphic as abstract fields. With the notations of the previous section $(E_A,\e_A)$ is the reflex field of $\MT(T_{A,\Q})$. We have an isomorphism of polarized $\Z$-HS $\phi \colon T_X \lr T_A$. We also look at the corresponding isomorphisms
$$
 \phi_\Q \colon T_{X,\Q} \lr T_{A,\Q} \ \ \ \text{and} \ \ \ \phi_f \colon T_{X,\A_f} \lr T_{A,\A_f}
$$
induced by the actions of $\pi,\iota,\gamma$ and $\beta$ on Betti cohomology with $\Q$ coefficients and on \'etale cohomology.
We have that $\phi_\Q = \phi\otimes_\Z \Q$ and $\phi_f = \phi \otimes_\Z \A_f$.

The morphism $\phi_\Q$ gives an isomorphism 
$$
\phi_\Q^{\rm ad} \colon E_X = \End_{\rm HS}(T_{X,\Q}) \lr \End_{\rm HS}(T_{A,\Q}) = E_A. 
$$
We have further the two inclusions $\e_X \colon E_X \lr \End_\C(H^{2,0}(X)) \cong \C$ and $\e_A \colon E_A \lr \End_\C(H^{2,0}(A)) \cong \C$. The map $\phi_\Q$ is defined as a composition of algebraic morphisms and hence we have a commutative diagram
\begin{equation}
\xymatrix{
E_X \ar[d]_{\e_X} \ar[rrrr]^{\phi_\Q^{\rm ad}} & & & & E_A \ar[d]_{\e_A} \\
\End_\C(H^{2,0}(X)) \ar[rrrr]^{(\pi^{* -1}_{DR} \circ \iota^*_{DR}\circ \gamma^{* -1}_{DR}\circ \beta^*_{DR})^{\rm ad}} & & & & \End_\C(H^{2,0}(A)).
}
\end{equation}
In other words $\phi_\Q^{\rm ad}$ gives an isomorphism of the CM-types $(E_X,\e_X)$ and $(E_A,\e_A)$. Therefore, with these identifications, $\phi_\Q$ commutes with the action of $E$ on the two vector spaces $T_{X,\Q}$ and $T_{A,\Q}$ via the isomorphisms $e_X^{-1} \colon E \lr E_X$ and $\e_A^{-1} \colon E \lr E_A$. Similarly $\phi_f$ is an $\A_{E,f}$-equivariant isomorphism via these actions. Further, via the isomorphism $\phi_\Q$, one can identify the cocharacters $\mu_X$ and $\mu_A$ and thus also the morphism $r_X$ and $r_A$. Taking all these remarks in to account we see that for any id\`ele $s \in \A^\times_{E}$ the following diagram
\begin{equation}\label{CommEtCoh1}
\xymatrix{
 T_{A,\A_f} \ar[r]^{r_A(s)} \ar[d]^{\phi_f} & T_{A,\A_f} \ar[d]^{\phi_f} \\
 T_{X,\A_f} \ar[r]^{r_X(s)} & T_{X,\A_f}
 }
\end{equation}
is commutative.

Let $\sigma$ be an element of $\Aut(\C/E)$. Making a base change $\Spec(\sigma) \colon \Spec(\C) \lr \Spec(\C)$ of Diagram \eqref{ShIConstruction} we obtain a diagram
\begin{equation}
\xymatrix{
& \tilde A^\sigma \ar[dl]_{\pi^\sigma} \ar[dr]^{\iota^\sigma} & & \tilde X^\sigma \ar[dl]_{\gamma^\sigma} \ar[dr]^{\beta^\sigma} & \\
A^\sigma & & Y^\sigma & & X^\sigma.
}
\end{equation}
Denote by $\phi^\sigma \colon T_{X^\sigma} \lr T_{A^\sigma}$ the isomorphism of polarized $\Z$-HS, defined as in Theorem \ref{ShIEK3-1}, by $(\pi^\sigma)^{* -1}\circ (\iota^\sigma)^* \circ (\gamma^\sigma)^{* -1}\circ (\beta^\sigma)^*$. Consider the isomorphisms induced on Betti cohomology with $\Q$ coefficients and \'etale cohomology
$$
 \phi_\Q^\sigma \colon T_{X^\sigma,\Q} \lr T_{A^\sigma,\Q} \ \ \ \text{and} \ \ \ \phi_f^\sigma \colon T_{X^\sigma,\A_f} \lr T_{A^\sigma,\A_f}
$$
These isomorphism are defined by algebraic morphisms and hence just as above we conclude that $\phi_\Q^\sigma$ defines an isomorphism of the reflex fields $(E_{X^\sigma},\e_{X^\sigma})$ and $(E_{A^\sigma},\e_{A^\sigma})$.
\newline
\newline
{\it Proof of Lemma \ref{RefFTwist}.} For an element $\sigma \in \Aut(\C/E)$ the reflex fields of $\Hg(T_{A,\Q})$ and $\Hg(T_{A^\sigma,\Q})$ are $E$. Hence we have an isomorphism $e_{A^\sigma} \circ (\phi^\sigma_\Q)^{\rm ad} \colon E_{X^\sigma} \lr E$ and therefore the reflex field of $\Hg(X^\sigma)$ is $E \subset \C$.
\qed

Before giving the proof of Theorem \ref{MainThCMEK3} we shall make a final remark. The map $\phi_f$ is defined as a composition of the algebraic maps $\pi, \iota, \gamma, \beta$ and their inverses acting on \'etale cohomology. The isomorphism $\phi^\sigma_f$ is defined in the same way using the conjugates $\pi^\sigma, \iota^\sigma, \gamma^\sigma, \beta^\sigma$. Hence we see that we have the following commutative diagram of \'etale transcendental spaces:
\begin{equation}\label{CommEtCoh2}
\xymatrix{
T_{X,\A_f} \ar[d]_{\phi_f} \ar[r]^{\sigma_{X,f}} & T_{X^\sigma, \A_f} \ar[d]_{\phi^\sigma_f} \\
T_{A,\A_f} \ar[r]^{\sigma_{A,f}} & T_{A^\sigma, \A_f}.
}
\end{equation}
{\it Proof of Theorem \ref{MainThCMEK3}.} Let $s\in \A^\times_E$ be an id\`ele such that $\text{art}_E(s) = \sigma|_{E^{\rm ab}}$. By Proposition \ref{CMPrAS} we have an isogeny 
$$
 \eta \colon A^\sigma \lr A
$$
such that for the induced isomorphism $\eta_f^* \colon T_{A,\A_f} \lr T_{A^\sigma, \A_f}$ on
\'etale cohomology with $\A_f$-coefficients we have $\eta_f^*(r_A(s)t) = \sigma_{A,f}(t)$ for every $t \in T_{A,\A_f}$.

Using the isomorphisms of $\Q$-HS $\phi_\Q$ and $\phi^\sigma_\Q$ we obtain an isomorphism of $\Q$-HS 
$$
 \eta_X \colon T_{X,\Q} \lr T_{X^\sigma,\Q}
$$ 
defined as $\eta_X = \phi^\sigma_\Q \circ  \eta_\Q^* \circ \phi_\Q^{-1}$. In other words we define $\eta_X$ by completing the diagram
\begin{equation}
\xymatrix{
T_{X,\Q} \ar[d]^{\phi_\Q} \ar@{-->}[r]^{\eta_X} & T_{X^\sigma,\Q} \ar[d]^{\phi^\sigma_\Q} \\
T_{A,\Q} \ar[r]^{\eta_\Q^*} & T_{A^\sigma, \Q}
}
\end{equation}
where $\eta_\Q^*$ is the $\Q$-HS morphism induced by $\eta$ on Betti cohomology. Note that from the remarks made above the isomorphism $\eta_X$ is $E_X$-equivariant. We get an $\A_{E,f}$-linear isomorphism $\eta_{X,f} \colon T_{X,\A_f} \lr T_{X^\sigma,\A_f}$ by tensoring $\eta_X$ with $\A_f$. By the commutativity Diagrams \eqref{CommEtCoh1} and \eqref{CommEtCoh2}, and using the fact that $\phi_f = \phi_\Q \otimes \A_f$ and $\phi^\sigma_f = \phi_\Q^\sigma \otimes \A_f$ we see that the following diagram is commutative:
\begin{equation}
\xymatrix{
T_{A,\A_f} \ar[dd]^{\phi_{f}} \ar[dr]^{r_A(s)} \ar[drrr]^{\sigma_{A,f}} & & & \\
 & T_{A,\A_f} \ar[dd]^{\phi_{f}} \ar[rr]^{\eta_f^*} & & T_{A^\sigma, \A_f} \ar[dd]^{\phi^\sigma_{f}} \\
T_{X,\A_f} \ar[dr]^{r_X(s)} \ar[drrr]^{\sigma_{X,f}} & & & \\
 & T_{X,\A_f} \ar[rr]^{\eta_{X,f}} & & T_{X^\sigma,\A_f}.
}
\end{equation}
Hence we have a $\Q$-HS isomorphism $\eta_X \colon T_{X,\Q} \lr T_{X^\sigma,\Q}$ such that $\eta_{X,f}(r_X(s)t) = \sigma_{X,f}(t)$ for every $t \in T_{X,\A_f}$.\qed
\begin{rem}
Note that the morphisms $\eta_{X}$ and $\eta_{X,f}$ are induced by a cycle in $X^\sigma \times X$. Indeed, if $\Gamma_\eta \subset A^\sigma \times A$ is the graph of $\eta$, then the isomorphisms $\eta_X$ and $\eta_{X,f}$ are given by the cycle
$$
 Z = \bigl((\beta^\sigma,\beta) \circ (\gamma^\sigma, \gamma)^{-1} \circ (\iota^\sigma, \iota) \circ (\pi^\sigma,\pi)^{-1}\bigr)(\Gamma_\eta) \subset X^\sigma \times X
$$
as in \S 3 of \cite{SK-StConj}.
\end{rem}
\subsection{Some Special Points on $Sh_\Kg(G,\Omega^\pm)_\C$}
Let $d \in \N$ and let $\Kg \subset \SO(V_{2d})(\hat \Z)$ be a subgroup of finite index such that $\Kg \subset \Kg_n$ for some $n \geq 3$. In order to carry out our strategy for
proving that the morphism $j_{d,\Kg,\C} \colon \Fk_{2d,\Kg,\C} \lr Sh_\Kg(G,\Omega^\pm)_\C$ is defined over $\Q$ we need to find enough special points on
$Sh_\Kg(G,\Omega^\pm)_\C$ for which we can control the Galois action.
\begin{prp}
Let $E \subset \C$ be a quadratic imaginary field. Then the set of special points $[x,a]_\Kg \in Sh_\Kg(G,\Omega^\pm)_\C(\C)$ with reflex field $E$ is dense for the Zariski topology.
\end{prp}
\begin{proof}
Let $C$ be an elliptic curve over $\C$ with CM by $E$ and consider the product abelian surface $A = C \times C$. Let $P$ be a point of infinite order in $C$ and consider the divisor
$$
D_1 = P\times C + C\times P
$$
on $A$. As $D_1 = pr_1^*P + pr_2^*P$, where $pr_i \colon A \lr C$ is the projection morphism onto the $i$-th factor, we have that it is an ample divisor on $A$. Its self-intersection number $(D_1,D_1)_A$ is $2$. 

Let $X$ be the Kummer surface associated to $A$. Then $X$ is an exceptional K3 surface and the reflex field of $\MT(X)$ is exactly $E$. Let $\pi \colon \tilde A \lr A$ be the blowing-up of $A[2]$ and $\iota \colon \tilde A \lr X$ be the morphism of degree 2. Then the line bundle
$$
 \L := \O_X\bigr(\iota_*(\pi^*(D_1))\bigr)
$$ 
defines a quasi-polarization on $X$ and one easily computes that $(\L,\L)_X = 2$. Hence $\L$ is primitive. 

Let $P^2_B(X,\Z(1))$ be the primitive Betti cohomology group with respect to $c_1(\L)$. Fix an isometry $a \colon P^2_B(X,\Z(1)) \lr L_2$. Then we have a point $x := a\circ h_{X} \circ a^{-1}$ in $\Omega^\pm$. The Mumford-Tate group of the $\Q$-HS $x$ induced on $V_2$ is $a^{\rm ad}(\MT(X))$ and hence its reflex field is $E$. By the strong approximation theorem the orbit $G(\Q)\cdot x$ is dense in $\Omega^\pm$ hence the set of points $\{[x,a]_\Kg | \ a \in G(\A_f)\}$ is dense in 
$$
 Sh_\Kg(G,\Omega^\pm)_\C = G(\Q)\backslash \Omega^\pm \times G(\A_f)/\Kg.
$$
\end{proof}
\begin{rem}
Similar density results appear in various papers on the Torelli theorem for K3 surfaces. The difference with our situation is that in those papers one mainly works with the full period domain of dimension 20. We also mention Lemma 7.1.2 in \cite{A-MHC} which almost gives the result we need.
\end{rem}
\begin{cor}\label{densityExcK3}
Let $d \in \N$ and let $E \subset \C$ be a quadratic imaginary field. The set of points $x \in \Fk_{2d,\Kg,\C}$ corresponding to exceptional K3 surfaces $X$ of CM-type $(E_X,\e_X)$ such that $\e_X(E_X) = E$ is dense for the Zariski topology in $\Fk_{2d,\Kg,\C}$.
\end{cor}
\begin{proof}
We have an \'etale morphism $j_{d,\Kg,\C} \colon \Fk_{2d,\Kg,\C} \lr Sh_\Kg(G,\Omega^\pm)_\C$. According to the preceding proposition the set of points $[x,a]_\Kg \in Sh_\Kg(G,\Omega^\pm)_\C(\C)$ with reflex field $E$ is dense in $Sh_\Kg(G,\Omega^\pm)_\C$. Therefore the preimage of this set under $j_{d,\Kg,\C}$ in $\Fk_{2d,\Kg,\C}$ is also dense. It consists exactly of the exceptional K3 surfaces $X$ of CM-type $(E_X,\e_X)$ such that $\e_X(E_X) = E$, with a polarization of degree $2d$ and a level $\Kg$-structure.
\end{proof}
\subsection{Complex Multiplication for K3 Surfaces}\label{CMK3SSect}
We will prove here that the field of definition of the morphism $j_{d,\Kg,\C}$ is $\Q$. To do that we will use the density result for exceptional polarized K3 surfaces and Theorem \ref{MainThCMEK3} which establishes a relation between the Galois action on such a surface and its periods.

\begin{thm}\label{CMK3}
Let $d \in \N$ and let $\Kg \subset \SO(V_{2d})(\hat \Z)$ be a subgroup of finite index such that $\Kg \subset \Kg_n$ for some $n \geq 3$. Then the morphism $j_{d,\Kg, \C}
\colon \Fk_{2d,\Kg,\C} \lr Sh_\Kg(G,\Omega^\pm)_{\C}$ is defined over $\Q$. In other words one has an \'etale morphism 
\begin{displaymath}
 j_{d,\Kg} \colon \Fk_{2d,\Kg,\Q} \lr Sh_\Kg(G,\Omega^\pm)
\end{displaymath}
such that $j_{d,\Kg}  \otimes \C = j_{d,\Kg,\C}$.
\end{thm}
\begin{proof} We will divide the proof into three steps.
\newline
{\bf Step 1.} Let $x \in \Fk_{2d,\Kg,\C}$ be a point corresponding to an exceptional K3 surface with CM by $E$. We will show first that for every $\sigma \in \Aut(\C/E)$ we have $j_{d,\Kg,\C}(\sigma (x)) = \sigma(j_{d,\Kg,\C}(x))$.
 
Let $E \subset \C$ be a quadratic imaginary field, let $(X,\l, \a)$ be a
polarized exceptional K3 surface of CM-type $(E_X,\e_X)$ with a level $\Kg$-structure $\a$ such that $e_X(E_X) = E$. Then we have the triple 
$$
 \bigl((P^2_B(X,\Z(1)), h_X), \psi_X, \tilde \a\Kg\bigr)
$$
where $\tilde \a$ is a representative of the class $\a$.

Let $\tilde \a$ be a representative of the class $\a$ and let $a_X \colon P^2_B(X,\Z(1)) \lr L_{2d}$ be an isometry as in the definition of the morphism $j_{d,\Kg,\C}$ (see Step 1 of the proof of Proposition \ref{PerMorphism}). Via this isometry we have an inclusion of algebraic groups $a^{\rm ad}_X \colon \MT(X) \hookrightarrow G$. By the modular description of $Sh_\Kg(G,\Omega^\pm)_\C$ and the definition of $j_{d, \Kg, \C}$ we have
$$
 j_{d,\Kg, \C}\bigl((X,\l,\a)\bigr) = [a_X \circ h_X \circ a_X^{-1},\ a_X \circ \tilde \a]_\Kg.
$$

We have that $P^2_B(X,\Q(1)) = T_{X, \Q} \oplus A_X^\l$ where 
$$
 A_X^\l := c_1(\l)^\perp \subset A_{X,\Q} \subset H^2_B(X,\Q(1))
$$ 
as polarized $\Q$-HS. By definition the action of $E_X$ on $A_X^\l$ is trivial. The same decomposition holds for \'etale cohomology $P^2_{\rm et}(X,\A_f (1)) = T_{X, \A_f} \oplus A_{X,\A_f}^\l$ where 
$$
 A_{X,\A_f}^\l := c_1(\l)^\perp \subset A_{X,\A_f} \subset H^2_{\rm et}(X,\A_f(1)).
$$
By the comparison theorem between Betti and \'etale cohomology these two decompositions are compatible with tensoring with $\A_f$.

Let $\sigma \in \Aut(\C/E)$ and consider the conjugate
$\sigma(X,\l,\a) = (X^\sigma,\l^\sigma,\a^\sigma)$ defined by the base change $\Spec(\sigma)
\colon \Spec(\C) \lr \Spec(\C)$. The surface $X^\sigma$ is also exceptional and for it we have similar decompositions of $P^2_B(X^\sigma, \Q(1))$ and $P^2_{\rm et}(X^\sigma, \A_f(1))$. The base change morphism $\Spec(\sigma)$ induces a morphism $\sigma_{\Pic} \colon \Pic(X) \lr \Pic(X^\sigma)$, preserving the intersection forms on both spaces and sending $\l$ to $\l^\sigma$. Hence we obtain an isomorphism of polarized $\Q$-HS $\sigma_{\Pic, \Q} \colon A_{X,\Q}^\l \lr A_{X^\sigma, \Q}^\l$ and an isomorphism $\sigma_{\Pic, f} \colon \colon A_{X,\A_f}^\l \lr A_{X^\sigma, \A_f}^\l$ such that $\sigma_{\Pic,\Q} \otimes \A_f = \sigma_{\Pic, f}$. Note that by its very definition $\sigma_{\Pic, f}$ is nothing else but $\sigma_f$ restricted to $A_{X, A_f}^\l$.

By Theorem \ref{MainThCMEK3} there exists an isomorphism of polarized $\Q$-HS $\eta_X \colon T_{X,\Q} \lr T_{X^\sigma, \Q}$ such that $\eta^{-1}_X \circ \sigma_f(t) = r_X(s)(t)$ for every $t\in T_{X,\A_f}$. Define the isomorphism of polarized $\Q$-HS
$$
 \eta = \eta_X \oplus \sigma_{\Pic, \Q} \colon P^2_B(X, \Q(1)) \lr P^2_B(X^\sigma, \Q(1)).
$$
Then we obtain an isomorphism of primitive \'etale cohomology
$$
\eta_f := \eta \otimes \A_f \colon P^2_{\rm et}(X,\A_f(1)) \lr P^2_{\rm et}(X^\sigma,\A_f(1))
$$
for which, by the remarks made above, we have $\eta^{-1} \circ \sigma_f(t) = r_X(s)(t)$ for every $t \in P^2_{\rm et}(X,\A_f(1))$.

Consider the isometry 
$$
 a_{X^\sigma} = a_X\circ \eta^{-1} \colon P^2_B(X^\sigma, \Q(1)) \lr V_{2d}. 
$$ 
We have that $\tilde \a^\sigma = \sigma_f \circ \tilde \a$ and we will see that $a_{X^\sigma} \circ \tilde \a^\sigma \in G(\A_f)$ i.e., that we can use the marking $a_{X^\sigma}$ to compute the periods of $(X^\sigma, \l^\sigma, \a^\sigma)$. We compute
\begin{equation}\label{rxcomp}
\begin{aligned}
a_{X^\sigma} \circ \tilde \a^\sigma &= a_X \circ \eta^{-1} \circ \tilde \a = a_X \circ r_X(s) \circ \tilde \a \\ 
&= a_X \circ a_X^{-1} \circ r_{(a_X\circ h_X\circ a_{X}^{-1})}(s) \circ a_X \circ \tilde \a \\
&= r_{(a_X\circ h_X\circ a_{X}^{-1})}(s) \circ a_X \circ \tilde \a
\end{aligned}
\end{equation} 
and hence $a_{X^\sigma} \circ \tilde \a^\sigma$ belongs to $G(\A_f)$. Here the morphism 
$$
 r_{(a_X\circ h_X\circ a_{X}^{-1})} \colon \A_{E,f}^\times \lr a_X^{\rm ad}(\MT(X))(\A_f)
$$ 
is the homomorphism associated to the special pair $(a_X^{\rm ad}(\MT(X)),a_X \circ h_X \circ a_X^{-1})$ of $(G,\Omega^\pm)$ as in Definition \ref{rhm}. From the modular description of $Sh_\Kg(G,\Omega^\pm)_\C(\C)$ given in Proposition \ref{Modinterpr1} we see that
$$
 j_{d,\Kg, \C}\bigl((X^\sigma,\l^\sigma,\a^\sigma)\bigr) =
[a_{X^\sigma}\circ h_{X^\sigma}\circ a^{-1}_{X^\sigma},\ a_{X^\sigma}\circ \tilde \a^\sigma]_\Kg.
$$
Then using \eqref{rxcomp} we compute
\begin{align*}
j_{d,\Kg, \C}\bigl(\sigma(X,\l,\a)\bigr)  &= [a_{X^\sigma}\circ h_{X^\sigma}\circ a^{-1}_{X^\sigma},\ a_{X^\sigma}\circ \tilde \a^\sigma]_\Kg \\
&= [a_X\circ (\eta^{-1} \circ h_{X^\sigma} \circ \eta) \circ a^{-1}_X,\ a_X \circ \a^{-1} \circ \sigma_f \circ \tilde \a]_\Kg \\
&= [a_X \circ h_X \circ a^{-1}_X,\ r_{(a_X\circ h_X\circ a_{X}^{-1})}(s) \circ a_X \circ \tilde \a]_\Kg\\
&= \sigma\bigl(j_{d,\Kg, \C}((X,\l,\a))\bigr).
\end{align*}
Hence the action of $\Aut(\C/E)$ on the point $(X,\l,\a)$ commutes with $j_{d,\Kg, \C}$.
\newline
\newline
{\bf Step 2.}
For a fixed quadratic field $E \subset \C$ the set of polarized exceptional K3 surfaces $X$ of CM-type $(E_X,\e_X)$ with $\e_X(E_X) = E$ is Zariski dense in $\Fk_{d,\Kg,\C}$ (see Corollary \ref{densityExcK3}). According to Step 1 the action of $\Aut(\C/E)$ commutes with $j_{d, \Kg, \C}$ on that set. Hence it commutes with $j_{d, \Kg, \C}$ and by Proposition 13.1 in \cite{Mil-IShV} we conclude that $j_{d, \Kg, \C}$ is defined over $E$.
\newline
\newline
{\bf Step 3.}
Choose two quadratic imaginary fields $E_1 \subset \C$ and $E_2 \subset \C$ such that $E_1 \cap E_2 = \Q$. By the previous step we know that $j_{d ,\Kg, \C}$ is defined over $E_1$ and $E_2$. Hence it is defined over their intersection $\Q$ which is the reflex field of $Sh(G,\Omega^\pm)_\C$.
\end{proof}
As a corollary of the preceding theorem one can obtain an analogue of the main theorem for complex multiplication for abelian varieties (Theorem 11.2 in \cite{Mil-IShV}) for K3 surfaces with CM. Let $X$ be a K3 surface over $\C$ with CM by $E = \e_X(E_X)$. We denote by $r_X \colon \A_E^\times \lr \MT(X)(\A_f)$ the morphism associated to the pair $(\MT(X),h_X)$.
\begin{cor}[Complex multiplication for K3 surfaces]\label{MCMThmK3S1}
Let $X$ be a primitively polarized K3 surface over $\C$ of degree $2d$ with CM by a field $E$. Recall that we consider $E$ embedded in to $\C$ via $\e_X$. For every $\sigma \in \Aut(\C/E)$ and an id\`ele $s \in \A^\times_E$ such that ${\rm art}_E(s) = \sigma|_{E^{\rm ab}}$ there is an isomorphism of polarized $\Q$-HS
$$
\eta \colon P^2_B(X,\Q(1)) \lr P^2_B(X^\sigma, \Q(1))
$$
such that for $\eta_f = \eta \otimes \A_f \colon P^2_{\rm et}(X, \A_f(1)) \lr P^2_{\rm et}(X^\sigma, \A_f(1))$ we have $\eta_f(r_X(s)t) = \sigma_f(t)$ for every $t \in P^2_{\rm et}(X, \A_f(1))$.
\end{cor}
\begin{proof}
Let $\l$ be the polarization on $X$. We can introduce a level $3$-structure $\a$ on $(X,\l)$ so that $(X,\l,\a) \in \Fk_{2d,3,\C}(\C)$. Let 
$$
 a_X \colon P^2_B(X,\Z(1)) \lr L_{2d}
$$ 
and 
$$
 a_{X^\sigma} \colon P^2_B(X^\sigma,\Z(1)) \lr L_{2d}
$$ 
be two markings as in the construction of the morphism $j_{d,\Kg_3,\C}$ (cf. Step 1 of the proof of Proposition \ref{PerMorphism}). Then we have that
$$
 j_{d,\Kg_3, \C}\bigl((X,\l,\a)\bigr) = [a_X \circ h_X \circ a_X^{-1},\ a_X \circ \tilde \a]_{\Kg_3}
$$
and 
$$
 j_{d,\Kg_3, \C}\bigl((X,\l,\a)\bigr) = [a_{X^\sigma} \circ h_{X^\sigma} \circ a_{X^\sigma}^{-1},\  a_{X^\sigma} \circ \tilde \a^\sigma]_{\Kg_3}.
$$
Using Theorem \ref{CMK3} and the definition of a canonical model (cf. Section \ref{CanModSect}) we see that
\begin{equation}\label{cmk3general}
\bigl(q\cdot (a_X \circ h_X \circ a_X^{-1}), q a_X \circ r_X(s) \circ \tilde \a\bigr) = 
(a_{X^\sigma} \circ h_{X^\sigma} \circ a_{X^\sigma}^{-1},\ a_{X^\sigma} \circ \sigma_f \circ \tilde \a)
\end{equation} 
for some $q \in G(\Q)$. Hence comparing the first terms in \eqref{cmk3general} we see that
$$
 \eta:= a_{X^\sigma}^{-1} \circ q \circ a_X \colon \colon P^2_B(X,\Q(1)) \lr P^2_B(X^\sigma, \Q(1))
$$
defines an isometry of $\Q$-HS. Form the equality between the second terms we see that $q a_X \circ r_X(s) = a_{X^\sigma} \circ \sigma_f$ i.e., that $\eta_f \circ r_X(s) = \sigma_f$.
\end{proof}
Before stating our final result we will point out a difference between the approach to the theory of complex multiplication for abelian varieties given, for instance, in Chapters 10, 11 and 12 of \cite{Mil-IShV} and the one for K3 surfaces given in this chapter. In the case of abelian varieties one first proves an analogue of Corollary \ref{MCMThmK3S1} and then derives an analogue of Theorem \ref{CMK3} from it (cf. \S 4 in \cite{D-ShV1}). Here we do the opposite as we do not see a way to prove directly Corollary \ref{MCMThmK3S1}. The reason is the following: For a K3 surface $X$ with CM by $E_X$ and an automorphism $\sigma \in \Aut(\C/E_X)$ one has little control over the transcendental lattice $T_{X^\sigma}$ of $X^\sigma$, unless $X$ is exceptional.
\begin{rem}
The statement in Corollary \ref{MCMThmK3S1} can be given in a form not including any polarizations. With the notations as above for a K3 surface with CM by $E$ one simply gets a $\Q$-HS isometry $\a \colon H^2_B(X, \Q(1)) \lr H^2_B(X^\sigma, \Q(1))$ such that for the induced isomorphism $\eta_f \colon H^2_{\rm et}(X, \A_f(1)) \lr H^2_{\rm et}(X^\sigma, \A_f(1))$ on \'etale cohomology one has that $\eta_f(r_X(s)t) = \sigma_f(t)$ for any $t \in H^2_{\rm et}(X,\A_f(1))$.
\end{rem}
\begin{cor}\label{CMK3FieldofDef}
Every complex K3 surface with CM can be defined over a number field which is an abelian extension of its Hodge endomorphism field.
\end{cor}
\begin{proof}
Let $(X,\l)$ be a polarized K3 surface of degree $2d$ with CM by $E_X$ and choose a full level $3$-structure $\a$ on $(X,\l)$. We have an open embedding 
$$
 j_{d,\Kg^{\rm full}_3} \colon \Fk^{\rm full}_{2d,3,\Q} \hookrightarrow Sh_{\Kg_3^{\rm full}}(G,\Omega^\pm)
$$ 
of schemes over $\Q$. The point $(X,\l,\a)$ maps to a special point in $Sh_{\Kg_3^{\rm full}}(G,\Omega^\pm)_\C$ which according to \S 3.15 in \cite{D-ShV1} can be defined over an abelian extension of its reflex field $e_X(E_X)$. Therefore $X$ can be defined over an abelian extension of $E_X$. Note that one can give a description of the corresponding extension in terms of the reciprocity law as explained in {\it loc. cit.}.
\end{proof}
\begin{rem}
Let us mention that a similar result can be found in the literature. Shioda and Inose (Theorem 6 in \cite{Sh-I}) prove that any exceptional K3 surface can be defined over a number field. As we shall see below one can actually give such a number field explicitly and conclude that it is an abelian extension of the Hodge endomorphism field of the exceptional K3 surface. Piatetskij-Shapiro and Shafarevich prove, using the Torelli theorem for K3 surfaces, that every K3 surface with CM can be defined over a number field. We refer to Theorem 4 in \cite{PSS-ArthK3}
\end{rem}
\begin{exa}
Let $X$ be an exceptional K3 surface with CM by $E$. Let $C_1$ and $C_2$ be the two elliptic curves from Theorem \ref{ShIEk32} and denote by $j_1$ and $j_2$ their $j$-invariants. By the theory of complex multiplication for elliptic curves we know that $K_1=E(j_1,C_1[2])$ and $K_2 = E(j_2,C_2[2])$ are abelian extensions of $E$. We can see that all morphisms and surfaces involved in the construction described in Theorems \ref{ShIEK3-1} and \ref{ShIEk32} are defined over the composite $K_1K_2$. Hence $X$ is defined over $K_1K_2$ which is an abelian extension of $E$.
\end{exa}
\subsection{Final Comments}\label{FinalComm}
{\bf (A) Complex multiplication.} Corollaries \ref{MCMThmK3S1} and \ref{CMK3FieldofDef} are analogues to two of the main theorems of the theory of complex multiplication for abelian varieties (Chapter I, Corollary 5.5 in \cite{Mil-CanModMixShV}). Another important result of that theory is that every abelian variety with CM defined over a number field $K$ has potentially good reduction at every prime ideal of $K$. We wonder if a similar result holds for K3 surfaces with CM.
\begin{que}
Let $K$ be a number field and suppose given a K3 surface with CM over $K$. Does $X$ have potentially good reduction at every prime ideal of $K$?
\end{que}
One could follow the line of thoughts of the proof of Proposition 10.5 in \cite{Mil-IShV}. In this way we can see that for a prime $\mathfrak p$ of $K$ the inertia action of $I_{\mathfrak p}$ on $H^2_{\rm et}(X_{\bar \Q},\Q_l)$ factorizes through a finite group. To finish ``the proof'' we need a N\'eron-Ogg-Shafarevich-type criterion for potentially good reduction of K3 surfaces. To our knowledge, in general, this is an open problem. Such a criterion exists for discrete valuation rings of characteristic zero. This follows form the degeneration result of Kulikov (Theorem II and Theorem 2.7 in \cite{Kul-K3Deg}, and \cite{PP-K3Deg}).
\newline
\newline
{\bf (B) The period morphism.} One knows that the period morphisms used in \cite{Ast-K3} and \cite{Fri-Torelli} are dominant. Further, the complement of their images are divisors. We will show here that the same holds for $j_{d,\Kg,\C}$.

Recall that we have a decomposition
\begin{displaymath}
Sh_{\Kg}(G,\Omega^\pm)_\C = \coprod_{[g] \in \mathcal C} \Gamma_{[g]} \backslash \Omega^+,
\end{displaymath}
where $\mathcal C := G(\Q)_+\backslash G(\A_f)/\Kg$ and $\Gamma_{[g]}= G(\Q)_+ \cap g \Kg g^{-1}$ for some representative $g$ of $[g] \in \mathcal C$. Let $X$ denote the geometric connected component of the canonical model $Sh_\Kg(G,\Omega^\pm)$ corresponding to $\Gamma_{[1]} \backslash \Omega^+$. Denote by $E_X \subset \C$ its field of definition. It is an abelian extension of $\Q$ and one can see that
\begin{equation}\label{GalDecompShV}
Sh_\Kg(G,\Omega^\pm) = \coprod_{\sigma \in \Gal(E_X/\Q)} X^\sigma.
\end{equation}
According to Propositions 7 and 8 in Expos\'e XIII in \cite{Ast-K3} we have that the intersection $j_{d,\Kg,\C}(\Fk_{2d,\Kg,\C})\cap X$ is dense in $X$ and its complement in $X$ is a divisor. The points in the complement correspond to quasi-polarized K3 surfaces. Hence using Theorem \ref{CMK3} and \eqref{GalDecompShV} we conclude that $j_{d,\Kg,\C}(\Fk_{2d,\Kg,\C})$ is dense in $Sh_\Kg(G,\Omega^\pm)_\C$ and its complement is a divisor.
\bibliographystyle{hplain}

\end{document}